\documentclass[11pt
]{article}
\usepackage{amsmath,amsthm,amssymb,color,epic,eepic, amsfonts,graphicx,mathrsfs,
cite
}

\usepackage{tikz}

\newlength{\minitwocolumn}
\newcommand{\Z}{{\Bbb Z}} 
\newcommand{\R}{{\Bbb R}} 
\newcommand{\C}{{\Bbb C}} 
\newcommand{\N}{{\Bbb N}} 
\newcommand{\FF}{{\Bbb F}} 
\newcommand{\T}{{\Bbb T}} 
\newcommand{\PP}{{\Bbb P}} 
\newcommand{\F}{{\mathcal F}}

\newcommand{\ta}{\tilde{a}}

\newcommand{\cD}{{\mathcal D}}

\newcommand{\cU}{{\mathcal U}}

\newcommand{\cI}{{\mathcal I}}

\newcommand{\cN}{{\mathcal N}}

\newcommand{\cR}{{\mathcal R}}

\newcommand{\cE}{{\mathcal E}}
\newcommand{\cM}{{\mathcal M}}
\newcommand{\cO}{{\mathcal O}}
\newcommand{\cQ}{{\mathcal Q}}

\newcommand{\cV}{{\mathcal V}}
\newcommand{\cW}{{\mathcal W}}

\newcommand{\bR}{{\overline{R}}}

\newcommand{\la}{\lambda}

\newcommand{\al}{\alpha}

\newcommand{\ep}{\epsilon}
\newcommand{\vep}{\varepsilon}

\newcommand{\tW}{\widetilde{W}}

\newcommand{\tit}{{\tilde{t}}}

\newcommand{\bep}{\bar{\epsilon}}

\newcommand{\bI}{\bar{I}}

\newcommand{\s}{{\sigma}}

\newcommand{\hf}{\widehat{f}}

\newcommand{\hV}{{\cV}}

\newcommand{\hrho}{
{\rho}}

\newcommand{\nn}{{\nonumber}}
\newcommand{\bea}{\begin{eqnarray}}
\newcommand{\ena}{\end{eqnarray}}
\newcommand{\be}{\begin{eqnarray*}}
\newcommand{\en}{\end{eqnarray*}}

\newcommand{\beit}{\begin{itemize}}
\newcommand{\enit}{\end{itemize}}

\newcommand{\lb}[1]{\label{#1}}

\newcommand{\ds}[1]{{\displaystyle #1 }}

\newcommand{\End}{\mathrm{ End}}
\newcommand{\rank}{\mathrm{ rank}}
\newcommand{\rk}{\mathrm{ rk}}

\newcommand{\id}{\mathrm{ id}}

\newcommand{\E}{\mathrm{E}}

\newcommand{\Stab}{\mathrm{Stab}}

\newcommand{\gC}{\mathfrak{C}}

\newcommand{\gS}{\mathfrak S}

\newcommand{\Wg}{{\cal W}\hspace{-0.1mm}\bar{\g}}
\newcommand{\Wppg}{{\cal W}_{p,p^*}\hspace{-0.1mm}\bar{\g}}

\newcommand{\bra}[1]{\langle #1 |}        
\newcommand{\ket}[1]{{| #1 \rangle}}      

\def\infq4p#1{{(#1;q^4,p)_\infty}}

\newcommand{\tot}{\, \widetilde{\otimes}\, }

\newcommand{\tr}{{\rm tr}}

\newcommand{\mmatrix}[1]{\begin{matrix} #1 \end{matrix}}

\font\teneufm=eufm10
\font\seveneufm=eufm7
\font\fiveeufm=eufm5
\newfam\eufmfam
\textfont\eufmfam=\teneufm
\scriptfont\eufmfam=\seveneufm
\scriptscriptfont\eufmfam=\fiveeufm

\let\goth\frak
\newcommand{\slth}{\widehat{\goth{sl}}_2}
\newcommand{\slt}{\goth{sl}_2}

\newcommand{\slnh}{\widehat{\goth{sl}}_N}
\newcommand{\sln}{\goth{sl}_N}
\newcommand{\g}{\goth{g}}

\newcommand{\Aqp}{{\mathcal A}_{q,p}}

\newcommand{\Bqla}{{{\mathcal B}_{q,\lambda}}}

\newcommand{\gl}{{\goth{gl}}}

\newcommand{\glnh}{\widehat{\goth{gl}}_N}

\newcommand{\h}{\goth{h}}

\font\fourteeneufm=eufm10 scaled\magstep2    
   
\newcommand{\gbig}{\mbox{\fourteeneufm g}} 
\makeatletter
\@addtoreset{equation}{section}
\makeatother

\begin{document}

\vspace{-1cm}
\begin{center}
{\bf\Large  Elliptic Quantum Groups\\[7mm] }
{\large  Hitoshi Konno }\\[6mm]
{\it  Department of Mathematics, Tokyo University of Marine Science and 
Technology, \\Etchujima, Koto, Tokyo 135-8533, Japan\\
       hkonno0@kaiyodai.ac.jp} 
\\[7mm]
\end{center}

\begin{abstract}
We expose the elliptic quantum groups 
in the Drinfeld realization associated with both the affine Lie algebra $\g$ and the toroidal algebra $\g_{tor}$. There the level-0 and level $\not=0$ representations appear in a unified way so that one can define the vertex operators as intertwining operators of them. 
The vertex operators are key for many applications such as  
a derivation of the elliptic weight functions, integral solutions of the (elliptic) $q$-KZ equation and a formulation of algebraic analysis of the elliptic solvable lattice models.  Identifying the elliptic weight functions with the elliptic stable envelopes we make a correspondence between  the level-0 representation of the elliptic quantum group 
  and  the equivariant elliptic cohomology. We also emphasize a characterization of the elliptic quantum groups as $q$-deformations of the $W$-algebras. 
\end{abstract}


\medskip
\noindent
{\bf Key points} \ :\ 
\begin{itemize}
\item The elliptic quantum groups $U_{q,p}(\g)$ and $U_{q,\kappa,p}(\g_{tor})$ are  formulated in terms of the Drinfeld generators  in a uniform way\footnote{One may unify them as ``quiver elliptic quantum groups'', where $U_{q,p}(\g)$ is 
a finite Dynkin quiver type, whereas $U_{q,\kappa,p}(\g_{tor})$  is an affine Dynkin quiver type. 
} and equipped with $H$-Hopf algebroid structure as a coalgebra structure.  
\item $U_{q,p}(\g)$ is characterized  both as a $p$-deformation of the quantum affine algebra $U_q(\g)$ and as a $q$-deformation of the $W$-algebra of the coset type $\g\oplus\g \supset \g_{\rm diag}$.
In the same way, $U_{q,\kappa,p}(\g_{tor})$  is characterized  both as a $p$-deformation of the quantum toroidal algebra $U_{q,\kappa}(\g_{tor})$ and as a $(q,\kappa)$-deformation of affine quiver $W$-algebra. 
\item The vertex operators defined as intertwiners of the tensor product of the level-0 and $\not=0$ modules are key objects for many applications.
\item The Gelfand-Tsetlin basis for the level-0 representation of $U_{q,p}(\slnh)$ is constructed explicitly by using the elliptic weight function as the change of basis matrix 
from the standard basis. 
 \item 
 The elliptic weight function of type $\sln$ is identified with the elliptic stable envelope for the equivariant elliptic cohomology of the cotangent bundle to the partial flag variety $\E_T(T^*fl_\la)$.  This allows us to make a dictionary between  $U_{q,p}(\slnh)$ and  $\E_T(T^*fl_\la)$.  
\end{itemize}

\section{Introduction}
Quantum groups are algebraic systems  associated with solutions of the Yang-Baxter equation 
(YBE). 
They are in general associative algebras classified by 
 finite dimensional simple Lie algebras $\bar{\g}$, affine Lie algebras $\g$ or toroidal algebras $\g_{tor}$, and equipped with certain coalgebra structure.  Typical examples are Yangians $Y(\bar{\g})$, Yangian doubles $\cD Y(\bar{\g})$, quantum groups $U_q(\bar{\g})$,  affine quantum groups $U_q(\g)$ and 
quantum toroidal algebras $U_{q,\kappa}(\g_{tor})$. See for example \cite{CPBook,Molev,GKV}.

Similarly elliptic quantum groups (EQG) are  quantum groups associated with elliptic  solutions of the YBE
 classified by affine Lie algebras $\g$ or  toroidal algebras $\g_{tor}$. 
In elliptic setting, it is known  at least for the cases classified by $\g$ that there are two types of YBEs, the vertex type and the face type 
(Fig.1 and Fig.2).  
Accordingly, there are two types of elliptic quantum groups, the vertex type and the face type\cite{JKOStg}. 
Both of them are dynamical quantum groups, whose universal $R$ matrices\cite{JKOStg} satisfy the dynamical YBE (DYBE)\cite{Felder}\footnote{ Obviously the face type YBE is dynamical\cite{Felder}. The vertex type is also dynamical, because the elliptic nome $p$ is a dynamical parameter getting a shift by a central element $q^{-2c}$ and yields the new elliptic nome $p^*=pq^{-2c}$\cite{JKOStg}.}.  
\vspace{5mm}

\begin{figure}[htbp]\lb{fig:vertexYBE}
\begin{center}
 \begin{tikzpicture}[scale=0.7]

   \draw [->,thick](2,-1)node[anchor=west]{$\vep_3$} --(-2,-1) node[anchor=east]{$\vep_3'$} node[anchor=south]{$\quad u_3$};
  
   \draw [->,thick](-1,0) node[anchor=south]{$\vep_1$}--(-1,-2)
;
   \draw [->, thick,](1,0) node[anchor=south]{$\vep_2$} --(1,-2)
   ;
   \draw [->,thick](-1,-2)
--(1,-3) node[anchor=north]{$\vep_1'$} node[anchor=west]{$u_1$};
   \draw [->, thick,](1,-2)
--(-1,-3) node[anchor=north]{$\vep_2'$} node[anchor=east]{$u_2$};

   \node (A) at (3.5,-1.5){$=$};

   \draw [->,thick](9,-2)node[anchor=west]{$\vep_3$} --(5,-2) node[anchor=east]{$\vep_3'$} node[anchor=south]{$\quad u_3$};
  
   \draw [->,thick](6,-1)
--(6,-3)
   node[anchor=north]{$\vep_2'$} node[anchor=east]{$u_2$}
;
   \draw [->, thick,](8,-1)
  --(8,-3)
    node[anchor=north]{$\vep_1'$}  node[anchor=west]{$u_1$}
   ;
   \draw [->,thick](6,0)
   node[anchor=south]{$\vep_1$}
--(8,-1)
;
   \draw [->, thick,](8,0)
    node[anchor=south]{$\vep_2$}
--(6,-1) ;

\end{tikzpicture}

{\footnotesize
Figure 1: The vertex type Yang-Baxter equation
}

\end{center}
\end{figure}

\begin{figure}[htbp]\lb{fig:faceYBE}
\begin{center}

\hspace{0.5cm}
\begin{tikzpicture}[scale=1.1]

\draw[thick] (-1,0) rectangle (0,1) ;
\draw[thick] (-0.8,1) -- (-1,0.8) ;

\draw[thick] (0,0) rectangle (1,1) ;
\draw[thick] (0.2,1) -- (0,0.8) ;

\draw[thick,rotate=-135] (0,0) rectangle (1,1) ;
\draw[thick,rotate=-135] (0.8,0) -- (1,0.2)   ;

\draw[dashed,thick] (-1,0) -- (-0.7,-0.7) ;
\draw[dashed,thick] (1,0) -- (0.7,-0.7) ;

\node[anchor=south] (A) at (-1,1){$a$};
\node[anchor=south] (B) at (0,1){$b$};
\node[anchor=south] (C) at (1,1){$c$};
\node[anchor=east] (F) at (-1,0){$f$};
\node[anchor=west] (D) at (1,0){$d$};
\node[anchor=north] (E) at (0,-1.4){$e$};
\node (G) at (0,0){$\bullet$};

\node (U1) at (-0.5,0.5){$u_1$};
\node (U2) at (0.5,0.5){$u_2$};
\node (U3) at (0,-0.7){{\footnotesize $\ u_1-u_2$}};

\node (EQ) at (2.2,0){$=$};

\draw[thick] (3.5,-1) rectangle (4.5,0) ;
\draw[thick] (3.7,0) -- (3.5,-0.2) ;

\draw[thick] (4.5,-1) rectangle (5.5,0) ;
\draw[thick] (4.7,0) -- (4.5,-0.2) ;

\draw[thick,rotate around={135:(4.5,0)}] (4.5,-1) rectangle (5.5,0) ;
\draw[thick,rotate  around={135:(4.5,0)}] (5.3,0) -- (5.5,-0.2)   ;

\draw[dashed,thick] (3.5,0) -- (3.8,0.7) ;
\draw[dashed,thick] (5.5,0) -- (5.2,0.7) ;

\node[anchor=east] (A1) at (3.5,0){$a$};
\node[anchor=south] (B1) at (4.5,1.4){$b$};
\node[anchor=west] (C1) at (5.5,0){$c$};
\node[anchor=north] (F1) at (3.5,-1){$f$};
\node[anchor=north] (D1) at (5.5,-1){$d$};
\node[anchor=north] (E1) at (4.5,-1){$e$};
\node (G1) at (4.5,0){$\bullet$};

\node (U12) at (4,-0.5){$u_2$};
\node (U22) at (5,-0.5){$u_1$};
\node (U32) at (4.5,0.7){{\footnotesize $\ u_1-u_2$}};

\end{tikzpicture}

{\footnotesize
Figure 2: The face type Yang-Baxter equation
}
\end{center}
\end{figure}

There are several formulations of the elliptic quantum groups. They can be classified by their generators and coalgebra structures. See Table \ref{tab:1} : $\Aqp(\slnh)$ and $\Bqla(\g)$
in terms of the Chevalley generators, $U_{q,p}(\g)$
in terms of the Drinfeld generators 
 and $E_{q,p}(\glnh)$
  in terms of the $L$ operators. 
Only $\Aqp(\slnh)$ is the vertex type.  The others are the face type.
Coalgebra structures are  the quasi-Hopf algebra structure
for $\Aqp(\slnh)$, $\Bqla(\g)$, 
 and the Hopf algebroid structure
 for $E_{q,p}(\glnh)$, 
  $U_{q,p}(\g)$, respectively. 
For expositions of these formulations with historical remarks and arguments for their consistency, the reader can consult\cite{KonnoBook}.

\begin{table}[htbp]
\begin{center}
\begin{tabular}{|c|c|c|}
\hline
& {co-algebra structure} & {generators}\\ \hline
&&\\[-4mm]
$\mmatrix{\Aqp(\slnh)\ \mbox{(vertex type)}\\[3mm]
         \Bqla({\g})\ \mbox{(face type)}\\[3mm]
}$ & quasi-Hopf algebra & Chevalley\\ \hline
&&\\[-3mm]
$\mmatrix{E_{q,p}(\g)\ \mbox{(face type)}\\[3mm]
 }$ & Hopf Algebroid & 
$L$-operator
\\ \hline
&&\\[-3mm]
   $\mmatrix{U_{q,p}({\g})\ \mbox{(face type)}\\[3mm]
 }$ & Hopf Algebroid & Drinfeld
 \\ \hline
\end{tabular}

\end{center}
\caption{Three formulations of the elliptic quantum groups}
\label{tab:1} 
\end{table}

In this article, we expose the elliptic quantum group $U_{q,p}(\g)$ and its toroidal version  $U_{q,\kappa,p}(\g_{tor})$. They are generated by the Drinfeld generators, i.e. analogues of the loop generators of the affine Lie algebra $\g$\cite{Kac}. 
After giving definitions of $U_{q,p}(\g)$ and some typical representations of $U_{q,p}(\slnh)$, we introduce the vertex operators of $U_{q,p}(\slnh)$ and describe their applications : derivations of the elliptic weight functions and integral solutions of the elliptic $q$-KZ equation, construction of the level-0 action on the Gelfand-Tsetlin basis,  algebraic analysis of the elliptic solvable lattice models. We then describe a geometric interpretation of
 the level-0 representation of $U_{q,p}(\slnh)$
 in terms of the equivariant elliptic cohomology of the cotangent bundle to the partial flag variety $\E_T(T^*fl_\la)$. We also emphasize a characterization of $U_{q,p}(\g)$ as a $q$-deformation of the $W$-algebra of the coset type. In the final section we briefly describe a formulation of the elliptic quantum toroidal algebras, $U_{q,\kappa,p}(\g_{tor})$, $U_{q,t,p}(\gl_{1,tor})$ and expose their connections to the Macdonald theory, affine quiver $W$-algebra  and application to  instanton calculus of the super symmetric gauge theories.

\section{Elliptic Algebra $U_{q,p}(\gbig)$}
Through this paper, let $p,q$ be generic complex numbers satisfying $|p|, |q|<1$. 
Let $\bar{\g}$ be a simple Lie algebra over $\C$ and $\g=X^{(1)}_l$  the corresponding untwisted affine Lie algebra with the generalized Cartan 
matrix $A=(a_{ij})_{i,j \in I}$, $I=\{0\}\cup \bar{I},\ \bar{I}=\{1,\dots,l\}$, \ $\rank A=l$. 
We denote by $B=(b_{ij})_{i,j \in I}, \,b_{ij}=d_i a_{ij}$ the symmetrization of $A$ and 
 set $q_i=q^{d_i}$. We also use the notations for $n \in \mathbb{Z}$,
\begin{eqnarray*}
&&[n]_q=\frac{q^n-q^{-n}}{q-q^{-1}},\qquad [n]_i=\frac{q_i^n-q_i^{-n}}{q_i-q_i^{-1}}. 
\end{eqnarray*}
We fix a realization $(\h,\Pi,\Pi^\vee)$ of $A$, i.e. $\h$ is a $l+2$-dimensional $\mathbb{C}$-vector space, 
$\Pi=\{\alpha_0,\alpha_1, \dots, \alpha_l\} \subset \h^*$ a set of simple roots, and $\Pi^\vee=\{h_0,h_1,\dots,h_l\} \subset \h$ a set of simple coroots satisfying 
$\langle\alpha_j,h_i\rangle=a_{ij}\ (i,j\in I)$ for a canonical pairing $\langle\, ,\rangle : \h^*\times \h \to \C$\cite{Kac}.  Set also $\cQ=\sum_{i\in I}\Z \al_i$. We take $\{h_1,\dots,h_l, c, d\}$  as the basis of $\h$ and $\{ \bar{\Lambda}_1,\cdots, \bar{\Lambda}_l, \Lambda_0, \delta \}$ the dual basis satisfying
\be
&&
\langle\delta,d\rangle=1=\langle\Lambda_0,c\rangle,\quad 
\langle\bar{\Lambda}_i,h_j\rangle=\delta_{i,j}\lb{pairinghhs}, 
\en
with the other pairings being 0. 

Let $H_P$ be  a $\C$-vector space spanned by $ P_0, P_1,\cdots,P_l$
 and $H^{*}_P$ be its dual space spanned by $Q_0, Q_1,\cdots,Q_l$ with a 
 pairing $\langle Q_i,P_j\rangle =a_{ij}$.   For $\al=\sum_{i\in I}c_i\al_i\in \h^*$, we set 
 $P_\al=\sum_{i\in I}c_iP_i$, $Q_\al=\sum_{i\in I}c_iQ_i$. In particular $P_{\al_i}=P_i$, $Q_{\al_i}=Q_i$. Define also an analogue of the root lattice $\cR_Q=\sum_{i\in I}\Z Q_i$. 
 
Let us consider $H:=\h\oplus H_P$ and  $H^*:=\h^*\oplus H^*_P$ 
 with a pairing $\langle \ ,\ \rangle : H^*\times H \to \C$ by extending those on $\h^*\times \h$ and $H^*_P\times H_P$ with $\langle \h^* ,H_P \rangle=0=\langle H^*_P ,\h \rangle$. 
We denote by $\FF=\cM_{H^*}$ the field of meromorphic functions on $H^*$. We regard 
a meromorphic function $g(h,P)$ of $h\in \h, P\in H_P$ as an element in $\FF$ by 
$g(h,P)(\mu)=g(\langle \mu,h\rangle, \langle \mu,P\rangle)$ for $\mu\in H^*$. 
 
\subsection{Definition}
The elliptic  algebra $U_{q,p}(\g)$ is a topological algebra over  $\FF[[p]]$  
generated by the Drinfeld generators 
\be
\al_{i,m}, \quad e_{i,n}, \quad f_{i,n}, \quad K_i^{\pm}, \quad q^{\pm c/2}, \quad q^{d}\quad 
 (i\in \bI, \quad m \in \mathbb{Z}\backslash\{0\}, \quad n \in \mathbb{Z}). 
\en
In order to write down the defining  relations, it is convenient to introduce their  generating functions $e_i(z), f_i(z)$ and $\phi^\pm_i(z)$ called the elliptic currents. 
\be
&&e_i(z)=\sum_{n\in \Z}e_{i,n} z^{-n},\qquad f_i(z)=\sum_{n\in \Z}f_{i,n} z^{-n},\\
&&\phi^+_i(q^{-{c}/{2}}z)
=K^+_i \exp\left( -(q_i-q_i^{-1}) \sum_{n>0} \frac{ \al_{i, -n}}{1-p^n}z^{ n}\right)
 \exp\left( (q_i-q_i^{-1}) \sum_{n>0} \frac{p^n\al_{i, n}}{1-p^n}z^{ -n}\right),\\
&&\phi^-_i(q^{{c}/{2}}z)
=K^-_i \exp\left( -(q_i-q_i^{-1}) \sum_{n>0} \frac{ p^n\al_{i, -n}}{1-p^n}z^{ n}\right)
 \exp\left( (q_i-q_i^{-1}) \sum_{n>0} \frac{ \al_{i, n}}{1-p^n}z^{ -n}\right).
\en
We also set $p^*=pq^{-2c}$. The defining relations are given as follows. 
For $g(h,P) \in \FF$, 
\be
&& q^{\pm c/2}
\  :\hbox{ central },\\
&&
g(h,{P})e_j(z)=e_j(z)g(h+\langle \al_j,h\rangle, P-\langle Q_{j},P\rangle),\quad
g(h,{P})f_j(z)=f_j(z)g(h-\langle \al_j,h\rangle, P),
\\
&&[g(h,P), \al_{i,m}]=
0,\qquad [g(h, P), q^d]=0,\quad 
g(h,{P})K^{\pm}_j=K^{\pm}_jg(h, P-\langle Q_{j},P\rangle),
\lb{gKpm}
\\
&& [q^d, \al_{j,m}]=q^m\al_{j,m},\quad [q^d, x^\pm_j(z)]=x^\pm_j(q^{-1}z)
\lb{dedf}\\
&&
[K^\pm_i,K^\pm_j]=[K^\pm_i,K^\mp_j]=0=[K^\pm_i,\al_{j,m}], \label{qta2}\\ 
&& K^\pm_i e_j(z)(K^\pm_{i})^{-1}=q^{\mp\langle \al_j,h_i \rangle }e_j(z), \quad K^\pm_i f_j(z)(K^\pm_{i})^{-1}=q^{\pm\langle\al_j,h_i\rangle }f_j(z),\label{qta3}
\en
\be
&&[\al_{i,m},\al_{j,n}]=\delta_{m+n,0}\frac{[a_{ij} m]_i}{m}\frac{q^{cm}-q^{-cm}}{q_j-q_j^{-1}}\frac{1-p^m}{1-p^{*m}}
q^{-cm},\\
&&
[\al_{i,m},e_j(z)]=\frac{[a_{ij} m]_i}{m}\frac{1-p^m}{1-p^{*m}}q^{-cm}
z^m e_j(z),\\
&&
[\al_{i,m},f_j(z)]=-\frac{[a_{ij} m]_i}{m}  
z^m f_j(z),
\\
&&(
z-q^{b_{ij}}w)g_{ij}(
w/z;p^*)
e_i(z)e_j(w)
=(
q^{b_{ij}}z-w) g_{ij}(
z/w;p^*)
e_j(w)e_i(z),\nn\\
&&(
z-q^{-b_{ij}}w)g_{ij}(
w/z;p)^{-1}
 f_i(z)f_j(w)
=(
q^{-b_{ij}}z-w) g_{ij}(
z/w;p)^{-1}
f_j(w)f_i(z),\nn\\
&&[e_i(z),f_j(w)]=\frac{\delta_{i,j}}{q_i-q_i^{-1}}
\left(\delta\bigl(q^c {w}/{z}\bigr)\phi^{-}_{i}(q^{{c}/{2}} w)
-\delta\bigl(q^{-c} {w}/{z}\bigr)\phi^{+}_{i}(q^{{c}/{2}}z)
\right),\lb{xpxm}
\\[2mm]
&&\mbox{+ Serre relations. }
\en
Here  we set
\be
&&g_{ij}(z;s)=\exp\left(-\sum_{m>0}\frac{1}{m}\frac{q^{b_{ij}m}-q^{-b_{ij}m}}{1-s^m}(sz)^m\right)\ \in \ \C[[s]][[z]].\lb{deg:g}
\en
These relations are treated as formal Laurent series in the argument of the elliptic currents i.e. $z, w$ etc.. All the coefficients in $z,w$ etc. are well defined in the $p$-adic topology. 
 
It is  easy to find that in the limit $p\to 0$ the above relations except for those indicating the non-commutativity of $\FF$ and $e_i(z), f_i(z), K^\pm_j$ go to the defining relations of the quantum affine algebra $U_q(\g)$\cite{DrinfeldNew}. One also finds that the non-commutativity of $\FF$ can be realized as the one of 
$P_j$ and $Q_i$ by setting $[P_i,Q_j]=a_{i,j}\ (i,j\in I)$.  Hence by using the group algebra   $\C[\cR_Q]$ of  $\cR_Q$, i.e. $e^{Q_\al}, e^{Q_\beta},  e^{Q_\al}e^{Q_\beta}=e^{Q_\al+Q_\beta}, e^0=1\in \C[\cR_Q]$,  
one obtains the following isomorphism\cite{FKO}. For generic $p,q$, 
\be
&&U_{q,p}(\g)/pU_{q,p}(\g)\ \cong \ 
(U_{q}(\g)\otimes \FF[[p]])\sharp \C[\cR_Q].
\en
Here the smash product $\sharp$ expresses the non-commutativity between $\FF$ and $\C[\cR_Q]$.  

\medskip
\noindent
{\it Remark.}  For representations, on which $q^{\pm c/2}$ take complex values e.g. $q^{\pm k/2}$ (see Sec.\ref{secRep}), we treat $p$ and  $p^*=pq^{-2k}$ as generic complex numbers satisfying  $|p|<1$ and $|p^*|<1$. 
Then  one has 
\be
 &&g_{ij}(z;s)=\frac{(sq^{b_{ij}}z;s)_\infty}{(sq^{-b_{ij}}z;s)_\infty},\qquad (z;s)_\infty=\prod_{n=0}^\infty(1-zs^n)
\en
for $|sq^{\pm b_{ij}}z|<1$, $s=p, p^*$.
Hence in the sense of analytic continuation, one can rewrite the relations of $e_i(z)$ and $e_j(w)$ 
and of $f_i(z)$ and $f_j(w)$ as
\bea
&&z\theta_{p^*}(q^{b_{ij}}
w/z)e_i(z)e_j(w)
=-w
\theta_{p^*}(q^{b_{ij}}
z/w)e_j(w)e_i(z),\lb{xpxp}\\
&&z\theta_{p}(q^{-b_{ij}}
w/z)f_i(z)f_j(w)
=-w
\theta_{p}(q^{-b_{ij}}
z/w)f_j(w)f_i(z),\lb{xmxm}
\ena
respectively. Here $\theta_p(z)$ denotes the odd theta function given by 
\be
\theta_{p}(z)=(z;p)_\infty(p/z;p)_\infty(p;p)_\infty.
\en

The additive notations are also often used.  Introduce {$r, r^*=r-k \ \in \R_{>0}$} and set $p=q^{2r}, p^*=pq^{-2k}=q^{2r^*}$, 
\be
&&
E_i(z):=e_i(z) {z^{-\frac{P_{\al_i}-1}{r^*}}},\quad F_i(z):=f_i(z) {z^{\frac{(P+h)_{\al_i}-1}{r}}}.
\en
Then from \eqref{xpxp}, \eqref{xmxm} and the non-commutativity of $h, P$ with $e_i(z), f_i(z)$  one obtains  
{
\be
&&F_i(z_1)F_j(z_2)=\frac{[u_1-u_2-a_{ij}/2]}{[u_1-u_2+a_{ij}/2]}F_j(z_2)F_i(z_1),\\
&&E_i(z_1)E_j(z_2)=\frac{[u_1-u_2+a_{ij}/2]^*}{[u_1-u_2-a_{ij}/2]^*}E_j(z_2)E_i(z_1).
\en
}
Here we set $z_i=q^{2u_i}\ (i=1,2)$. The symbols $[u]$ and $[u]^*$ denote 
Jacobi's odd theta functions given by
\be
&&[u]=q^{\frac{u^2}{r}-u}\theta_p(q^{2u}), 
\quad [u]^*=q^{\frac{u^2}{r^*}-u}\theta_{p^*}(q^{2u}).
\en
They have the following quasi periodicity. 
\be
&&[u+r]=-[u], \quad [u+r\tau]=-e^{-\pi i\tau}e^{-2\pi i {u}/{r}}[u],  \\
&&[u+r^*]^*=-[u]^*, \quad [u+r^*\tau^*]^*=-e^{-\pi i\tau^*}e^{-2\pi i {u}/{r^*}}[u]^*, 
\lb{thetaquasiperiod}
\en
where $p=e^{-2\pi i/\tau },\ p^*=e^{-2\pi i/\tau^*}$.

\subsection{Coalgebra structure
 } 
In order to formulate a coalgebra structure over $\FF=\cM_{H^*}$, which does not commute with the Drinfeld generators, one needs to extend the Hopf algebra to the $H$-Hopf algebroid. For the elliptic algebra  $\cU=U_{q,p}(\g)$, it is roughly sketched as follows. 
For expositions of the background concepts of $H$-Hopf algebroid, the reader can consult\cite{EV99,KR,Konno09}.

Let $P+h=\sum_{i\in I}c_i(P_i+h_i)\in H,\ c_i\in \C$. 
The algebra $\cU$ is {a $H$-algebra} by 
\be
&&\cU=\bigoplus_{\al,\beta\in \h^*}\cU_{\al,\beta}\\
&&(\cU)_{\al\beta}=\left\{x\in U \left|\ q^{P+h}x q^{-(P+h)}=q^{\langle\al,P+h\rangle}x,\quad q^{P}x q^{-P}=q^{\langle Q_\beta,P\rangle}x\quad \forall P+h, P\in H\right.\right\}
\en
with the two moment maps $\mu_l, \mu_r : \FF \to \cU_{0,0}$ defined by 
\be
&&\mu_l(\hf)=f(h,P+h,p)\in \FF[[p]],\qquad \mu_r(\hf)=f(h,P,p^*)\in \FF[[p^*]]
\en
for $\hf=f(h,P,p^*)\in \FF[[p^*]]$. Here  $p^*=pq^{-2c}$ as before.

The tensor product $\cU {\widetilde{\otimes}}\,\cU$ is the $\h^*$-bigraded vector space with 
\be
 (\cU {\widetilde{\otimes}}\,\cU)_{\al\beta}=\bigoplus_{\gamma\in \h^*} (\cU_{\al\gamma}\otimes_{\cM_{H^*}}\cU_{\gamma\beta}),
\en
where $\otimes_{\cM_{H^*}}$ denotes the ordinary tensor product 
modulo the following 
relation.
\be
&&\mu_r(\widehat{f})a\tot b=a\tot \mu_l(\widehat{f})b\qquad a, b\in \cU.\ \lb{fsabafb}
\en

Let us regard $T_\al=e^{-Q_\al}\in \C[\cR_Q]$ as a shift operator 
\be
&&(T_\al \widehat{f})=f(h,P+\langle Q_\al,P\rangle,p). 
\en
Then $\cD=\{\widehat{f} e^{-Q_\al} \ |\  \widehat{f}\in \FF,  e^{-Q_\al}\in \C[\cR_Q]\}$ becomes the $H$-algebra  having the property 
\be
&&\cU\cong \cU\tot \cD\cong  \cD\tot \cU \lb{Diso} 
\en
by $a\cong a\tot T_{-\beta}\cong T_{-\al}\tot a$ for all $a\in \cU_{\al\beta}$.

Let $\pi_V, V=\oplus_{i}\C v_i$ be the vector representation of $\bar{\g}$ and consider the elliptic dynamical $R$-matrix $R^+(z,\Pi)\in \End(V\otimes V)$ of  type $\g$\cite{JKOStg,Konno06}. 

For example, the  $\slnh$  type  is given by 
\begin{eqnarray}
R^+(z,\Pi)&=&\hrho^+(z)\bar{R}(z,\Pi),\lb{ellR}
\\
\bar{R}(z,\Pi)&=&
\sum_{j=1}^{N}E_{j,j}\otimes E_{j,j}+
\sum_{1 \leq j_1< j_2 \leq N}
\left(b_{}(u,(P+h)_{j_1,j_2 })
E_{j_1,j_1}
\otimes E_{j_2,j_2}+
\bar{b}_{}(u)
E_{j_2,j_2}\otimes E_{j_1,j_1}
\right.
 \nonumber\\
&&\qquad 
\left.
+
c_{}(u,(P+h)_{j_1,j_2 })
E_{j_1,j_2}\otimes E_{j_2,j_1}+
\bar{c}_{}
(u,(P+h)_{j_1,j_2 })E_{j_2,j_1}\otimes E_{j_1,j_2}
\right),\lb{ellR}\nn
\end{eqnarray}
where $E_{i,j}v_{k}=\delta_{j,k}v_i$, $z=q^{2u}$, $\Pi_{j,k}=q^{2(P+h)_{j,k}}$, 
$(P+h)_{j,k}=\sum_{i=j}^{k-1}(P_i+h_i)
$, 
\be
&&{\rho}^+(z)=q^{-\frac{N-1}{N}}z^{\frac{N-1}{rN}}
\frac{\Gamma(z;p,q^{2N})\Gamma(q^{2N}z;p,q^{2N})}{\Gamma(q^2z;p,q^{2N})\Gamma(q^{2N-2}z;p,q^{2N})}, 
\\[1mm]
&&b(u,s)=
\frac{[s+1][s-1][u]}{[s]^2[u+1]},\qquad 
\bar{b}(u)=
\frac{[u]}{[u+1]},\lb{Relements}\\
&&c(u,s)=\frac{[1][s+u]}{[s][u+1]},\qquad \bar{c}(u,s)=\frac{[1][s-u]}{[s][u+1]}.
\en
The symbol $\Gamma(z;p,q^{2N})$ denotes the elliptic Gamma function defined by
\bea
&&\Gamma(z;p,q)=\frac{(pq/z;p,q)_\infty}{(z;p,q)_\infty},\qquad (z;p,q)_\infty=\prod_{m,n=0}^\infty(1-zp^mq^n), \quad |p|, |q|<1.  \lb{ellGamma}
\ena

Let ${L}^+(z)=\sum_{i,j}E_{i,j}L^+_{i,j}(z)\in \End(V)\tot\cU[[z,z^{-1}]]$  be the the $L$-operator satisfying the dynamical $RLL$-relation\cite{JKOStg}
\begin{eqnarray*}
R^{+(12)}(z_1/z_2,\Pi){L}^{+(1)}(z_1)
{L}^{+(2)}(z_2)=
{L}^{+(2)}(z_2)
{L}^{+(1)}(z_1)
R^{+*(12)}(z_1/z_2,\Pi^*).\label{thm:RLL}
\end{eqnarray*}
Here $R^{+*}(z,\Pi^*)$ denotes the same elliptic $R$-matrix $R^{+}(z,\Pi^*)$ 
with  $\Pi^*_{i,j}=q^{2P_{i,j}}$ except for replacing the elliptic nome $p$ by $p^*$. 
Define two $H$-algebra homomorphisms, 
the counit $\vep : \cU\to \cD$ and the comultiplication $\Delta : \cU\to \cU \widetilde{\otimes}\; \cU$,  
and the antihomomorphism $S:\cU\to \cU$ by
\be
&&\vep(L^+_{i,j}(z))=\delta_{i,j}{T}_{Q_{\ep_i} }\quad (n\in \Z),
\qquad \vep(e^Q)=e^Q,\lb{counitUqp}\\
&&\vep(\mu_l({\hf}))= \vep(\mu_r(\hf))=\widehat{f}T_0, \lb{counitf}\\
&&\Delta(L^+_{i,j}(z))=\sum_{k}L^+_{i,k}(z)\widetilde{\otimes}
L^+_{k,j}(z),\lb{coproUqp}\\
&&\Delta(e^{Q})=e^{Q}\tot e^{Q},\\
&&\Delta(\mu_l(\hf))=\mu_l(\hf)\widetilde{\otimes} 1,\quad \Delta(\mu_r(\hf))=1\widetilde{\otimes} \mu_r(\hf),\lb{coprof}\\
&&S(L^+_{ij}(z))=(L^+(z)^{-1})_{ij},\\
&&S(e^{Q})=e^{-Q}, \quad S(\mu_r(\hat{f}))=\mu_l(\hat{f}),\quad S(\mu_l(\hat{f}))=\mu_r(\hat{f}).
\en
Here $\al_i=\ep_i-\ep_{i+1}\ (i\in \bI)$. Then the set  $(U_{q,p}(\g),H,{\cM}_{H^*},\mu_l,\mu_r,\Delta,\vep,S)$ becomes a $H$-Hopf algebroid\cite{Konno09,Konno18,KO,KonnoBook}. 

An explicit realization of $L^+(z)$ in terms of the elliptic currents of $\cU$  was studied in\cite{JKOS,KK03,Konno18,KO}.  In general, the connection between $L^+(z)$ and the elliptic currents is  given as follows. 
Let $\{h^1,\cdots,h^l\}$ be the dual basis to $\{h_1,\cdots,h_l\}$ of $\bar{\h}$ and $(\pi_V, V)$ be the vector representation of $\bar{\g}$.  Define\cite{JKOStg} 
\be
&&L^-(z)=({\rm Ad}(q^{-2\theta_V(P)})\otimes \id)(q^{2T_V}L^+(zpq^{-c}),\\
&&\theta_V(P)=-\frac{1}{2}\sum_{j}\left(\pi_V(h_j)\pi_V(h^j)+2P_j\pi_V(h^j)\right),\\
&&T_V=\sum_{j}\pi_V(h_j)\otimes h^j. 
\en
One finds that $L^\pm(z)$ satisfy 
\bea
&&R^{\pm(12)}(z_1/z_2,\Pi){L}^{\pm(1)}(z_1)
{L}^{\pm(2)}(z_2)=
{L}^{\pm(2)}(z_2)
{L}^{\pm(1)}(z_1)
R^{\pm*(12)}(z_1/z_2,\Pi^*),\label{RLL1}\\
&&R^{+(12)}(q^{c}z_1/z_2,\Pi){L}^{+(1)}(z_1)
{L}^{-(2)}(z_2)=
{L}^{-(2)}(z_2)
{L}^{+(1)}(z_1)
R^{+*(12)}(q^{-c}z_1/z_2,\Pi^*).\label{RLL2}
\ena
Here $R^{\pm}(z,\Pi), R^{\pm*}(z,\Pi^*)$ and $L^+(z)$ are related to $R^{\pm}_{VV}(z,\la+h), R^{\pm}_{VV}(z,\la)$ and $L^+_V(z,\la)$ in \cite{JKOStg} by
\be
&&R^{\pm}(z,\Pi)=R^{\pm}_{VV}(z,\la+h),\qquad R^{\pm*}(z,\Pi^*)=R^{\pm}_{VV}(z,\la),\\
&&  L^+(z)=L^+_{V}(z,\la)e^{-\sum_iE_{i,i}Q_{\ep_i}}
\en
with $\la=(r-k+h^\vee)d+\sum_i(P_{\al'_i}+1)\bar{h}^i$, $\la+h=(r+h^\vee)d+\sum_i(P_{\al'_i}+h_{\al'_i}+1)\bar{h}^i$, $[P_{i,j},Q_{\ep_k}]=\delta_{i,k}-\delta_{j,k}$ and $\al'_i$ being the simple root of the Langlands dual Lie algebra of $\bar{\g}$\cite{Konno06}. 
%
Then consider the following Gauss decomposition of $L^\pm(z)$. 
\be
&&L^\pm(z)=F^\pm(z)K^\pm(z)E^\pm(z). 
\en
Here 
$K^\pm(z)$ are diagonal matrices with entries $K^\pm_j(z)$, whereas $F^\pm(z)$ (resp. $E^\pm(z)$) are upper (resp. lower) triangular matrices with entries $F^\pm_{i,j}(z)$ (resp. $E^\pm_{j,i}(z)$) $(i<j)$ and  all diagonal  entries 1. By combining the relations for these entries obtained from \eqref{RLL1}-\eqref{RLL2}, one finds the following identification with the elliptic currents $E_j(z), F_j(z)$ and  $\phi^\pm_j(z)$. 
\be
&&E_j(zq^{j-c/2})=\mu^*\left(E^+_{j+1,j}(zq^{c/2})-E^-_{j+1,j}(zq^{-c/2})\right),\\
&&F_j(zq^{j-c/2})=\mu\left(F^+_{j,j+1}(zq^{-c/2})-F^-_{j,j+1}(zq^{c/2})\right),\\
&&\phi^\pm_j(zq^{- c/2}q^j)=\kappa K^\pm_j(z)K^\pm_{j+1}(z)^{-1}.
\en
Here $\mu, \mu^*, \kappa$  are constants satisfying
\be
&&\mu\mu^*=-\kappa\frac{q}{q-q^{-1}}\frac{(p^*;p^*)_\infty}{\theta_{p^*}(q^2)}.
\en
See for example \cite{Konno18,KonnoBook}. 

Note also that to formulate a $H$-Hopf algebroid structure 
one may also use the Drinfeld comultiplication for the elliptic currents\cite{KOgl1,KO23} instead of the standard comultiplication for $L^+(z)$ given in the above. The resultant $H$-Hopf algebroid is completely different coalgebra structure from the above. 
These two coalgebra structures have their own appropriate applications. 
In particular, they define different vertex operators as intertwining operators of 
 $U_{q,p}(\g)$-modules, even if the modules are the same. 
The standard one is appropriate to discuss the problems exposed in the 
subsequent sections Sec.\ref{secVO}-\ref{secGR}, whereas the Drinfeld comultiplication is appropriate to discuss problems related to the deformed $W$-algebras. See  Sec.\ref{UqpWpps} and Sec.\ref{secEQTA}. One should also note that  for the quantum toroidal  algebras, whatever they are trigonometric or elliptic, 
 only the coalgebra structures associated with the Drinfeld comultiplication are available  at the moment, because of a lack of expressions of the (elliptic) $R$-matrices and the $L$-operators.

\section{$U_{q,p}(\gbig)$ as Deformed $W$-algebra}\lb{UqpWpps}
The elliptic algebra $U_{q,p}(\g)$ has another characterization as a $q$-deformation of the $W$-algebra $\Wg
$ of the coset type $(\g)_{r-h^\vee-k}\oplus (\g)_{k}\supset (\g_{\rm diag})_{r-h^\vee}$\cite{GKO,BoSc,CrRav}, or more precisely as an algebra of screening currents of the deformation of $\Wg$. Here we assume $\g$ is an untwisted affine Lie algebra and realized as 
a central extension of the loop algebra\cite{Kac}. The generating functions of $\g$ is called the currents. The symbol $(\g)_k$ means to consider $\g$ at the level $k$ representation. 
A key for understanding this characterization of $U_{q,p}(\g)$ is  the (generalized) Feigin-Fuchs (FF) construction of 
$\Wg
$  from $\g$\cite{FFuchs,CrRav}. 
There one starts from the level-$k$ currents of  $\g$
  in the  Lepowsky-Wilson (LW) realization\cite{LW} and the corresponding energy-momentum (EM) tensor of the Wess-Zumino-Witten (WZW) model, i.e. the generating function of the Virasoro algebra.  The LW realization of $\g$ is a realization of 
the currents in terms of the $\rk\; {\g}$ number of free bosons and the level-$k$ $Z$-algebra. 
The EM tensor  is obtained by the Sugawara construction in terms of the currents. Then the  FF construction deforms the EM tensor by introducing the so-called background charge term, which depends on a parameter $r$, and at the same time deforms the free boson part of  the currents of $\g$ by $r$ and the $Z$-algebra part by adding the dynamical parameters. The deformed EM  tensor becomes the generating function of the Virasoro generators of the coset $W$-algebra $\Wg$,  whereas the deformed currents of $\g$ becomes the so-called screening currents, i.e. the primary field of  conformal dimension one w.r.t. this new Virasoro generators. 
All the generating functions of $\Wg$, the Virasoro as well as the other generators of higher conformal dimensions, are characterized as the intersection of the kernel of the screening operators i.e. certain contour integral of the screening currents, on $U(\g)$\cite{FeFr}. 
In this sense, we refer the resultant  $r$-deformation (+ modification by the dynamical parameters) of the  affine Lie algebra $\g$ to the $W$-algebra $\Wg$. When one  consider a $q$-deformation of it, $\g$ and its 
 $r$-deformation $\Wg$  turn out to give $U_q(\g)$ and its $p$-deformation $U_{q,p}(\g)$ with $p=q^{2r}$, respectively.  
 
 For example, the currents of the level $k(\not=0)$ affine Lie algebra $\slth$ has the following  LW realization.  
 \be
&&e(z)=Z^+(z):\exp\left\{\sqrt{2k}\sum_{m\not=0}\frac{a_m}{km}z^{-m}\right\}:,\quad f(z)=Z^-(z):\exp\left\{-\sqrt{2k}\sum_{m\not=0}\frac{a_m}{km}z^{-m}\right\}:, 
\en
where $Z^\pm(z)$ are the generating functions of the level-$k$ $Z$-algebra, and $a_m\ (m\in \Z\backslash\{0\})$ generate the Heisenberg subalgebra of $\slth$ satisfying
\bea
&&[\sqrt{{2k}}a_m,\sqrt{{2k}}a_n]=2mk\delta_{m+n,0}. \lb{aman}
\ena
The generating function $\widetilde{\phi(z)}=\sum_{m\not=0}\frac{a_m}{m}z^{-m}$ is called free boson.  
The $Z$-algebra is realized in terms of the $\Z_k$-parafermions  $\Psi(z), \Psi^\dagger(z)$ 
  and the zero mode operators $\al$, $h$ satisfying $[h,\al]=2$ as follows.
\bea
&&Z^+(z)=\Psi(z)e^{\al}z^{h/k},\qquad  Z^-(z)=\Psi^\dagger(z)e^{-\al}z^{-h/k}.\lb{Zop}
\ena
The corresponding 2d conformal field theory (CFT) is the $\slth$ WZW model generated by the EM  tensor 
 \be
&&T(z)=\frac{1}{2}(\partial\widetilde{\phi(z)})^2+\mbox{(zero-modes term)}+T_{PF}(z).
\en
 Here 
$T_{PF}(z)$ denotes the EM tensor of the 
$\Z_k$-parafermions. 
The central charge of the models is $c_{WZW}=\frac{3k}{k+2}$. 

Then the  FF construction is the following procedure. 
Deform $T(z)$ 
by adding the background charge term to $T(z)$  
 \be
&&
T(z) \ \mapsto \ 
T_{FF}(z)=T(z)
+\sqrt{\frac{k}{2r(r-k)}}
\partial^2\widetilde{\phi(z)}.
\en
This makes $T_{FF}(z)$ the generating function of  the coset Virasoro algebra $(\slth)_{r-2-k}\oplus (\slth)_k\supset (\slth)_{r-2}$ whose 
 central charge is 
$c_{Vir}=\frac{3k}{k+2}\left(1-\frac{2(k+2)}{r(r-k)}\right)$. 
At the same time, deform the currents of $\slth$ by scaling the coefficients of the free boson and by 
adding the dynamical parameter $P$ and its conjugate  $Q$ satisfying $[P,Q]=2$ to the $Z$-algebra part :
\bea
e(z)&\mapsto& S^+(z)= Z^+(z)e^{-Q}z^{-\frac{P-1}{r-k}}:\exp\left\{\sqrt{\frac{2kr}{r-k}}\sum_{m\not=0}\frac{a_m}{km}z^{-m}\right\}:,\lb{FFe}\\
f(z)&\mapsto& S^-(z)=Z^-(z)z^{\frac{P+h-1}{r}}:\exp\left\{-\sqrt{\frac{2k(r-k)}{r}}\sum_{m\not=0}\frac{a_m}{km}z^{-m}\right\}:. \lb{FFf}
\ena
It then turns out $S^\pm(z)$ become the screening currents of the   Virasoro algebra generated by $T_{FF}(z)$.

Now let us consider a $q$-deformation of these constructions. It is known that the quantum affine algebra  $U_q(\g)$ has an analogue of the LW realiozation. See for example \cite{FKO}. 
In our case, the level-$k$ currents of $U_q(\slth)$ are realized in terms of the Heisenberg generators $\ta_m$ and the level-$k$ $q$-deformed $Z$ operators $Z^\pm_q(z)$.   The former   
satisfy 
 \bea
&&[\ta_m,\ta_n]=\frac{[2m]_q[km]_q}{m}q^{-k|m|}\delta_{m+n,0}.\lb{almaln}
\ena
This is a $q$-deformation of \eqref{aman}, i.e.  in the limit $q\to 1$ one recovers   \eqref{aman}. The operators $Z^\pm_q(z)$ are realized in the same was as \eqref{Zop} by replacing $\Psi(z), \Psi^\dagger(z)$ with appropriate $q$-parafermions\cite{Konno,FKO}.  
The level-$k$ currents of $U_q(\slth)$ are hence realized as 
\be
&&e_q(z)=Z_q^+(z):\exp\left\{\sum_{m\not=0}\frac{\ta_m}{[km]_q}z^{-m}\right\}:,\quad f_q(z)=Z_q^-(z):\exp\left\{-\sum_{m\not=0}\frac{\ta_m}{[km]_q}z^{-m}\right\}:. 
\en

Then a $q$-analogue of the FF construction \eqref{FFe}-\eqref{FFf} is obtained  by noting 
the following observations. 
\begin{itemize} 
\item In \eqref{FFe} and \eqref{FFf}, the $Z$-algebra part remains the same as in $e(z)$ and $f(z)$  except for getting a modification by the dynamical parameters 
\item The commutation relation of the Heisenberg generators $a_m$ with the scaled coefficients can be read as from $S^+(z)$,  
\bea
&&[\sqrt{\frac{2kr}{r-k}}a_m,\sqrt{\frac{2kr}{r-k}}a_n]=2mk\frac{r}{r-k}\delta_{m+n,0}, \lb{camcane}
\ena
and  from $S^-(z)$, 
\bea
&&[\sqrt{\frac{2k(r-k)}{r}}a_m,\sqrt{\frac{2k(r-k)}{r}}a_n]=2mk\frac{r-k}{r}\delta_{m+n,0}. \lb{camcanf}
\ena
\end{itemize}  
The first point indicates that  the  $q$-$Z$-algebra operators $Z^\pm_q(z)$  in  $e_q(z)$ and $f_q(z)$  should remain the same except for getting a modification by the dynamical parameters $P$ and $Q$  in a $q$-deformed FF process.  
The second point suggests to use the Heisenberg generators $\al_m$ or $\al'_m$ satisfying the following  $p$-deformed commutation relation from \eqref{almaln}. 
\bea
&&[{\al}_m, {\al}_n]=\frac{[2m]_q[km]_q}{m}\frac{1-p^m}{1-p^{*m}}q^{-km}\delta_{m+n,0},\lb{comalale}\\
&&[{\al}'_m, {\al}'_n]=\frac{[2m]_q[km]_q}{m}\frac{1-p^{*m}}{1-p^{m}}q^{km}\delta_{m+n,0}.\lb{comalalf}
\ena
Here we set $p=q^{2r}$, $p^*=pq^{-2k}=q^{2(r-k)}$, and  $\al_m$, $\al'_m$ are not independent
\be
&&{\al}'_m=\frac{1-p^{*m}}{1-p^m}q^{km}{\al}_m.
\en 
These are $q$-deformations of \eqref{camcane} and \eqref{camcanf}.  
Note that \eqref{comalale} is nothing but 
the relations for the Heisenberg subalgebra of $U_{q,p}(\slth)$ 
 at level $k$ given in Sec. 2.1.  
One hence obtains the following $q$-analogue of the FF construction. 
\be
e_q(z)&\mapsto& S^+_q(z)= Z_q^+(z)e^{-Q}z^{-\frac{P-1}{r-k}}:\exp\left\{
\sum_{m\not=0}\frac{{\al}_m}{[km]_q}z^{-m}\right\}:,\\
f_q(z)&\mapsto& S^-_q(z)=Z_q^-(z)z^{\frac{P+h-1}{r}}:\exp\left\{-
\sum_{m\not=0}\frac{{\al}'_m}{[km]_q}z^{-m}\right\}:.  
\en
Then the currents $S^+_q(z)$ and $S^-_q(z)$ turn out to give a realization of the level-$k$ elliptic currents $e(z)$ and $f(z)$ of $U_{q,p}(\slth)$, respectively\cite{Konno,FKO}.  Hence regarding $e(z)$ and $f(z)$ as screening currents, one can define a deformation of the $W$-algebra ${\cal W}\slt$ of the coset type at level $k$ as the kernel of the screening operators associated to either $e(z)$ or $f(z)$. 

In general, in the FF process  the $Z$-algebra structure of $\g$ remains the same as in the LW realization to realize the screening currents of $\Wg
$ of the coset type $(\g)_{r-h^\vee-k}\oplus (\g)_{k}\supset (\g_{\rm diag})_{r-h^\vee}$ except for a modification by adding the dynamical parameters. Correspondingly, the $q$-$Z$-algebra structure of $U_{q}(\g)$ remains the same to realize the elliptic currents $e_i(z), f_i(z)$ of $U_{q,p}(\g)$ as a $p$-deformation of $U_{q}(\g)$  except for a modification by the dynamical parameters\cite{Konno,FKO}. Then by using $e_i(z)$ and  $f_i(z)$ as screening currents, one defines a deformation of $\Wg$  
as the intersection of the kernel of the screening operators associated with either $e_i(z)$'s or $f_i(z)$'s. 
We denote the resultant deformed $W$-algebra by $\cW_{p,p^*}\hspace{-0.2mm}\bar{\g}$.   
\be
U_q(\g)&&\\[-2mm]  
{\footnotesize q\mbox{-deform.}}\ \rotatebox[origin=c]{45}{$\longrightarrow$}\hspace{1cm}&\rotatebox[origin=c]{-45}{$\longrightarrow$}&{\footnotesize {p=q^{2r}}\mbox{-deform.}
+\mbox{dynamical param.} }\\[-1mm]
&&\\[-3mm]
\g\hspace{2.3cm}&&\quad U_{q,p}(\g)\\
&&\\[-3mm]
{\footnotesize {r}\mbox{-deform.}+\mbox{dynamical param.} 
}\ \rotatebox[origin=c]{-45}{$\longrightarrow$}\hspace{1cm}&\rotatebox[origin=c]{45}{$\longrightarrow$}&{\footnotesize q\mbox{-deform.}}\\[-2mm]
\mbox{\footnotesize (the Feigin-Fuchs construction)}\qquad\qquad&&\\[-3mm]
\Wg
&&\\[-2mm]
&&\hspace{-1.8cm}\parallel\\[-3mm] 
&&
\hspace{-3.5cm}
(\g)_{r-h^\vee-k}\oplus (\g)_k\supset (\g)_{r-h^\vee}
\en
\begin{center}
{\footnotesize Figure 3 : Two chacterizations of $U_{q,p}(\g)$\qquad\qquad\qquad}
\end{center}

\medskip
Note that here we have obtained two elliptic nomes $p$ and $p^*$  naturally 
 in a $q$-deformation of the FF construction. 
 The parameters $p, p^*$ correspond to $q,t$ in Frenkel-Reshetikhin's 
formulation of the deformed $W$-algebra  $W_{q,t}(\bar{\g})$\cite{FrRe}\footnote{
Frenkel-Reshetikhin's  deformed $W$-algebras are the Hamiltonian reduction type, which are isomorphic to the coset type only for the simply laced $\bar{\g}$. For example, for 
the case $B_l$ the coset construction gives Fateev-Lukyanov's $W\hspace{-1mm}B_l$ algebra\cite{LuFa}, which is different from $W(B_l)$ obtained by the quantum Hamiltonian reduction of $B^{(1)}_l$. See for example \cite{BoSc}. }.  
It is remarkable that the appearance of two elliptic nomes is also consistent to the quasi-Hopf formulation of the face type EQG $\Bqla(\g)$\cite{JKOStg}, where $p$ is treated as a dynamical parameter and $p^*$ is nothing but a dynamical shift of $p$ by $q^{-2c}$.    

The generating function of $\cW_{p,p^*}\hspace{-0.2mm}\bar{\g}$ 
can also be constructed  explicitly in terms of the  vertex operators of $U_{q,p}(\g)$ w.r.t. the Drinfeld comultiplication. 
The same construction is valid also in the elliptic quantum toroidal algebras, which give a realization of the deformed affine quiver $W$-algebras. See Sec.\ref{secEQTA} for the case $\gl_{1,tor}$.  

Moreover in Sec.\ref{secVO}, we discuss the vertex operators of $U_{q,p}(\slnh)$
 derived as intertwining operators of $U_{q,p}(\slnh)$-modules  w.r.t. the standard comultiplication. 
Constructing correlation functions of them, one finds that they play the role of the vertex operators of the corresponding deformed $W$-algebra 
$\cW\hspace{-0.7mm}{}_{p,p^*}\hspace{-0.2mm}{\sln}$. 
See Sec.\ref{secSolqKZ} and \ref{secAA}.

\section{Representations}\lb{secRep}
A  representation of  $H$-algebra  is called a dynamical representation. 
We skip all formal descriptions of it and expose only the level-0 and the level-1 dynamical representations of $U_{q,p}(\slnh)$. Here one says that 
a $U_{q,p}(\g)$-module has level $k$, if $q^{c/2}$ acts  as the scalar $q^{k/2}$ on it.  
For detailed description of the dynamical representations, the reader can consult\cite{EV99,KR,Konno09}.

\subsection{The evaluation representation }\lb{vecrep}
Let  $\displaystyle{\hV=\oplus_{i=1}^{ N} \FF v_i\otimes 1}$ and set $\hV_z=\hV[z,z^{-1}]$. 
The following gives the level-0 dynamical action of $U_{q,p}(\slnh)$  
on $\hV_z$.   
\be
&&\pi_z(q^{c/2})=1,\quad \pi_z(d)=-z\frac{d}{dz},\\
&&\pi_z(\al_{j,m})=\frac{[m]_q}{m}(q^{j-N+1}z)^m(q^{-m}E_{j,j}-q^mE_{j+1,j+1}),\\
 && \pi_z(e_j(w))=\frac{(pq^2;p)_{\infty}}{(p;p)_{\infty}}E_{j,j+1}
     \delta\left(q^{j-N+1}{z}/{w}\right) 
    e^{-Q_{\al_j}}, \\
 && \pi_z(f_j(w))=\frac{(pq^{-2};p)_{\infty}}{(p;p)_{\infty}}E_{j+1,j}
     \delta\left(q^{j-N+1}{z}/{w}\right), \\
 && \pi_z(\phi_j^{+}(w))=q^{-\pi(h_j)}e^{-Q_{\al_j}}
     \frac{\theta_p(q^{-j+N-1+2\pi(h_j)} \frac{w}{z})}{\theta_p(q^{-j+N-1} {w}/{z})},\\
 && \pi_z(\phi_j^{-}(w))=q^{\pi(h_j)}e^{-Q_{\al_j}}
     \frac{\theta_p(q^{j-N+1-2\pi(h_j)} \frac{z}{w})}{\theta_p(q^{j-N+1} {z}/{w})},\qquad j=1,\cdots, N-1.
\en
Here $\pi(h_j)=E_{j,j}-E_{j+1,j+1}$, and the element $e^{Q}\in \C[\cR_Q]$ acts on $f(h,P)v\otimes 1$  by $e^{Q}\cdot(f(h,P)v\otimes 1)=f(h,P-\langle P,Q\rangle)v\otimes e^Q$. 
This  is called the evaluation 
representation associated with the vector representation with the evaluation parameter $z$.  
In particluar,  one has 
\be
&&\pi_z({L^+_{i,j}(w)})_{k,l}=R^+(w/z,\Pi^*)_{ik}^{jl}.
\en

\subsection{The level-1 highest weight representation}\lb{levelonerep}
Let $\Lambda_0$ and $\Lambda_a=\Lambda_0+\bar{\Lambda}_a \ (a=1,\cdots,N-1)$ be the fundamental weights of $\slnh$. For generic 
$\nu\in \h^*$, set 
\be
\hV(\Lambda_{a}+\nu,\nu)&=&\FF\otimes_\C (\F_{\al}\otimes e^{{\Lambda}_a}\C[{\cQ}])\otimes 
e^{Q_{{\nu}}}\C[{\cR}_Q],
\en 
where $\F_{\al}=\C[ \alpha_{j,-m}\ (j=1,\cdots,N-1,\ m\in \N_{>0})]$ denotes the Fock module of the Heisenberg algebra. 
We here consider the group algebra 
 $\C[{\cQ}]\otimes\C[\cR_Q]$ with the following central extension.
\be
&&e^{ \alpha_j}e^{ \alpha_k}=(-1)^{(\alpha_j, \alpha_k)} q^{\frac{1}{r}(\delta_{j,k+1}-\delta_{k,j-1})}e^{\alpha_k}e^{ \alpha_j},\\
&&
e^{Q_{\al_{j}}}e^{ Q_{\al_{k}}}=
q^{\left(\frac{1}{r}-\frac{1}{r^*}\right)(\delta_{j,k+1}-\delta_{j,k-1})} e^{ Q_{\al_{k}}}e^{Q_{\al_{j}}},\\
&&e^{Q_{\al_j}}e^{\al_k}=q^{\frac{1}{r}(\delta_{j,k+1}-\delta_{k,j-1})}e^{\al_k}e^{Q_{\al_j}}.
\en
Then the  space $\hV(\Lambda_{a}+\nu,\nu)$ is the level-1 irreducible highest weight $U_{q,p}(\slnh)$-module with the highest weight  $(\Lambda_a+\nu,\nu)$ w.r.t. the set $(P+h,P)$ by 
the  action
\begin{eqnarray*}
 && E_j(z) =\; :\exp \left\{ -\sum_{n \ne 0} \frac{1}{[n]_q} \alpha_{j,n}z^{-n} \right\}:
 e^{{\alpha}_j}e^{-Q_{\alpha_j}}z^{h_{\alpha_j}+1}(q^{N-j}z)^{-\frac{P_{\alpha_j}-1}{r^*}},  \\
 && F_j(z) =\; :\exp \left\{ \sum_{n \ne 0}\frac{1}{[n]_q} \alpha'_{j,n}z^{-n}\right\}: e^{-{\alpha}_j}z^{-h_{\alpha_j}+1}
(q^{N-j}z)^{\frac{(P+h)_{\alpha_j}-1}{r}}, \lb{level1Fj}\\
&&\phi^+_j(q^{-{1}/{2}}z)=q^{-h_j}e^{-Q_{\al_j}}:\exp\left\{(q-q^{-1})\sum_{m\not=0}\frac{\al_{j,m}}{1-p^m}p^mz^{-m}\right\}:,\\
&&\phi^-_j(q^{{1}/{2}}z)=q^{h_j}e^{-Q_{\al_j}}
:\exp\left\{(q-q^{-1})\sum_{m\not=0}\frac{\al_{j,m}}{1-p^m}z^{-m}\right\}:
\quad  (1\leq j\leq N-1).
\end{eqnarray*}
The highest weight vector is given by $1\otimes  e^{\Lambda_a}\otimes e^{Q_{{\nu}}}$. 

Note that one has the following natural decomposition.  
\be
&&\hV(\Lambda_{a}+\nu,\nu)=\bigoplus_{\xi, \eta\in \cQ}
\F_{a,\nu}(\xi,\eta),\qquad 
\F_{a,\nu}(\xi,\eta)=\FF\otimes_\C(\F_{\al}\otimes e^{{\Lambda}_a+{\xi}})\otimes  e^{Q_{{\nu}+\eta}}.\lb{repsp}
\en
It turns out that $\F_{a,\nu}(\xi,\eta)$ is isomorphic to the Verma module of the $W$-algebra 
$\cW\sln$ with the central charge $(N-1)\left(1-\frac{N(N+1)}{r(r-1)}\right)$ and the highest weight\\ $\frac{1}{r(r-1)}\bigl|r(\nu+\eta+\rho)-(r-1)(\Lambda_a+\nu+\xi+\eta+\rho)\bigr|^2$\cite{FKO}.

\section{Vertex Operators }\lb{secVO}
The  vertex operators defined for the infinite dimensional quantum groups such as $U_q(\g)$,  $U_{q,p}(\g)$ and $U_{q,\kappa,p}(\g_{tor})$ have proven to be key objects in many applications\cite{IFrRe,JM,KKW,CKSW2,KOgl1}. 
Here we expose the vertex operators of $U_{q,p}(\slnh)$\cite{JKOS,KK03}, which are used 
in a derivation of the elliptic weight functions (Sec.\ref{secWF}), integral solutions of the elliptic $q$-KZ equation and the vertex function (Sec.\ref{secSolqKZ}) and also in the algebraic analysis of the elliptic lattice models (Sec.\ref{secAA}).

Let ${\cV}_z$ be as in Sec.\ref{vecrep} 
and ${\cV}(\Lambda_a+\nu,\nu)$ denote the irreducible 
level-$1$ highest weight $U_{q,p}(\slnh)$-module with the highest weight $(\Lambda_a+\nu,\nu)$ in 
Sec.\ref{levelonerep}. Let $\Lambda_{a+N}=\Lambda_a$.  
The (type I)  vertex operator $\Phi(z)$ is the intertwinner of 
the $U_{q,p}(\slnh)$-modules
\be
 {\Phi}(z) &:& {\cV}(\Lambda_a+\nu,\nu) \to  {\cV}_z \tot {\cV}(\Lambda_{a-1}+\nu,\nu)
\en
satisfying the intertwining relations 
\bea
 && \Delta(x) {\Phi}(z) ={\Phi}(z)x      
 \qquad \forall x \in U_{q,p}(\slnh).\label{irI}
\ena
The components of the vertex operator are defined by
\begin{equation*}
 {\Phi}(zq^{-1})u=\sum_{\mu=1}^{ N}v_{\mu} 
 \tot \Phi_{\mu} \left(z \right)u, \quad
 \qquad  \forall u \in {\cV}(\la,\nu).
\end{equation*}

One can solve the linear equation \eqref{irI} by using the representations in the last section 
uniquely up to an overall constant factor. The result is summarized as follows. 
\bea
\Phi_{\mu}(z)
&=&a_{\mu,N}\oint_{\T^{N-\mu}}\prod_{m=\mu}^{N-1}\frac{dt_m}{2\pi i t_m}\Phi_{N}(z)F_{N-1}(t_{N-1})F_{N-2}(t_{N-2})\cdots F_{\mu}(t_{\mu})
\varphi_{\mu}(z,t_{\mu},\cdots,t_{N-1};\Pi),\nn\\
 \Phi_{N}(z) &=& : \exp \left\{ \sum_{m \ne 0} (q^m-q^{-m}){\cE}_m^{'N} z^{-m} \right\}:
  e^{-{\bep}_N
} (- z)^{-h_{\bep_N}}z^{\frac{1}{r}(P+h)_{\bep_N}}, \label{typeIvo}
\ena
where we set $\bep_i=\ep_i-\sum_{j=1}^N\ep_j/N$, $\Pi=\{\Pi_{\mu,m}\ (m=\mu+1,\cdots,N)\}$,   
\be
&&
\T^{N-\mu}=\{t\in \C^{N-\mu}\ |\ |t_\mu|=\cdots=|t_{N-1}|=1 \}, 
\en
and 
\be
&&{\varphi_{\mu}(z,t_{\mu},\cdots,t_{N-1}; \Pi)=
\prod_{m=\mu}^{N-1}
\frac{[v_{m+1}-v_{m}+(P+h)_{\mu,m+1}-\frac{1}{2}][1]}
{[v_{m+1}-v_{m}+\frac{1}{2}][(P+h)_{\mu,m+1}]}}\lb{defvarphi}
\en
with $z=q^{2u}, t_m=q^{2v_m}$, $v_N=u$. We assume $|p|<|z|<1$. 
The symbol $\cE_m^{'j}$ denotes a linear combination of $\al'_{j,m}$ obtained  by solving 
\be
&&\al'_{j,m}=[m]_q^2(q-q^{-1})(\cE_m^{'j}-q^{-m}\cE_m^{'j+1}),\qquad \sum_{j=1}^Nq^{(j-1)m}\cE^{' j}_m=0.
\en
   

The most important property of  $\Phi_\mu(z)$ is the following commutation relation.
\begin{eqnarray}
 \Phi_{\mu_2}(z_2)\Phi_{\mu_1}(z_1) &=& \sum_{\mu_1',\mu_2'=1}^{N}R(z_1/z_2,\Pi)_{\mu_1 \mu_2}^{\mu_1' \mu_2'}\ 
 \Phi_{\mu_1'}(z_1)\Phi_{\mu_2'}(z_2), \label{typeIcr}
\end{eqnarray}
where 
\be
&& R(z,\Pi)=\mu(z)\bR(z,\Pi),\qquad 
\mu(z)=z^{-\frac{r-1}{r}\frac{N-1}{N}} \frac{\Gamma(pz;p,q^{2N})\Gamma(q^{2N}z;p,q^{2N})}{\Gamma(q^2z;p,q^{2N})\Gamma(pq^{2N-2}z;p,q^{2N})}.
\en
For example, this yields the transition property of the elliptic weight function (Sec.\ref{secWF}) as well as the $R$-matrix coefficients in the elliptic $q$-KZ connection (Sec.\ref{secSolqKZ}).

\section{Elliptic Weight Functions}\lb{secWF}
The weight functions are objects in the theory of 
($q$-) hypergeometric integrals. They play the role of a basis of the (twisted) de Rham cohomology\cite{AK,Matsuo,Mimachi,TV97}. The elliptic version of the weight function was  first introduced for the $\slt$ type in \cite{TV97}\footnote{Strictly speaking, the elliptic weight functions in \cite{TV97} were used as a pole subtraction matrix, in the terminology of \cite{AO}, in the $q$-hypergeometric integral solution of the $q$-KZ equation. They specifies the cycles of the integral. On the other hand, there the role of a basis of the co-cycle was  played by the trigonometric weight functions. }.  We here expose the elliptic weight functions of the $\sln$ type, which can be derived by considering a composition of the vertex operators given in the last section. In the next section, 
we show that they play the role in elliptic hypergeometric integral solutions of the elliptic $q$-KZ equation. Moreover, the elliptic weight functions are identified with the elliptic stable envelopes\cite{AO} associated with Nakajima quiver variety $X$. This makes clear a connection between representation theory of EQG and geometry of the elliptic cohomology (Sec.\ref{secGR}). 
  
Let  $\Phi_\mu(z)$ be  the vertex operators of $U_{q,p}(\slnh)$ in the last section  and consider their  composition
\be
&&\phi_{\mu_1\cdots\mu_n}(z_1,\cdots,z_n):=\Phi_{\mu_1}(z_1)\cdots\Phi_{\mu_n}(z_n)\ :\ 
\F_{a,\nu}(\xi,\eta) \ \to \ \F_{a',\nu}(\xi,\eta).
\en
Here $a'\equiv a-n$ mod $N$. 
%
Let $[1,n]=\{1,\cdots,n\}$ and 
define the index sets $I_l:=\{ i\in [1,n]\ |\ \mu_i=l\}$ $(l=1,\cdots,N)$. Set also $\la_l:=|I_l|$, 
$\la:=(\la_1,\cdots,\la_N)$.  Then $I=(I_1,\cdots,I_N)$ is a partition of $[1,n]$, i.e.  
 \be
I_1\cup \cdots \cup I_N=[1,n],\quad I_k\cap I_l=\emptyset\quad  \mbox{$(k\not=l)$}.
\en
Thus obtained partition $I$ is often denoted by $I_{\mu_1,\cdots\mu_n}$.
For $\la=(\la_1,\cdots,\la_N)\in \N^N$, let $\cI_\la$ be the set of all partitions  $I=(I_1,\cdots,I_N)$ satisfying $|I_l|=\la_l\ (l=1,\cdots,N)$.
 Set also $\la^{(l)}:=\la_1+\cdots+\la_l$,  $I^{(l)}:=I_1\cup\cdots \cup I_l$ and let $I^{(l)}=:\{i^{(l)}_1< \cdots<i^{(l)}_{\la^{(l)}}\}$. Note that $\cI_\la$ 
 specifies the coordinate flags of the partial flag variety $fl_\la$ : $0=V_0\subset V_1\subset \cdots \subset V_{N-1}\subset V_N=\C^n$ with $\dim V_l=\la^{(l)}$. 

Now let us substitute the realization of $\Phi_\mu(z)$ \eqref{typeIvo} to 
 $\phi_{\mu_1\cdots\mu_n}(z_1,\cdots,z_n)$. In each  $\Phi_\mu(z)$, we assign the integration variable $t^{(l)}_a$ to the elliptic current $F_l(\ast)$ $(l=1,\cdots,N-1)$ as its argument if $\Phi_\mu(z)$ is the $i^{(l)}_a$-th vertex operator i.e. $z=z_{i^{(l)}_a}$.  
 After rearranging the order of the elliptic currents $F_l$'s and $\Phi_N$'s and taking their normal ordering as specified in $\widetilde{\Phi}(t,z)$ in the below, one obtains the following expression for $I=I_{\mu_1\cdots\mu_n}$,  
\bea
&&\phi_{\mu_1\cdots\mu_n}(z_1,\cdots, z_n)=\oint_{\T^M
} \prod_{l=1}^{N-1}\prod_{a=1}^{\la^{(l)}}\frac{dt^{(l)}_a}{2\pi i t^{(l)}_a}
\  \widetilde{\Phi}(t,z)\tW_I
(t,z,\Pi),\qquad |p|<|z_1|,\cdots, |z_n|<1, \nn\\
&&\lb{nptop}
\ena
\be
&&\widetilde{\Phi}(t,z)= :\Phi_N(z_1)
\cdots 
\Phi_N(z_n):
:F_{N-1}(t^{(N-1)}_{1})\cdots F_{N-1}(t^{(N-1)}_{\la^{(N-1)}}):\cdots  :F_1(t_1^{(1)})\cdots F_{1}(t^{(1)}_{\la^{(1)}}):
\nn\\
&&
\qquad\qquad
\times \prod_{1\leq k<l\leq n}<\Phi_N(z_k)\Phi_N(z_l)>^{Sym}
\prod_{l=1}^{N-1}\prod_{1\leq a<b\leq \la^{(l)}}<F_l(t^{(l)}_a)F_l(t^{(l)}_b)>^{Sym},
\lb{Phitilde}\\
&&
\widetilde{W}_I(t,z,\Pi)= {\rm Sym}_{t^{(1)}}\cdots {\rm Sym}_{t^{(N-1)}}
\widetilde{U}_I(t,z,\Pi),\lb{def:Wtilde}
\en
\be
\widetilde{U}_I(t,z,\Pi)
&=&\prod_{l=1}^{N-1}\prod_{a=1}^{\la^{(l)}}\left(\frac{[v^{(l+1)}_b-v^{(l)}_a+(P+h)_{\mu_s,l+1}-C_{\mu_s,l+1}(s)][1]}{[v^{(l+1)}_b-v^{(l)}_a+1]\left[(P+h)_{\mu_s,l+1}-C_{\mu_s,l+1}(s)\right]}\right|_{i^{(N)}_s=i^{(l+1)}_b=i^{(l)}_a}\\[2mm]
&&
\qquad\qquad\qquad\times \prod_{b=1\atop i^{(l+1)}_b>i^{(l)}_a}^{\la^{(l+1)}}\frac{[v^{(l+1)}_b-v^{(l)}_a]}{[v^{(l+1)}_b-v^{(l)}_a+1]}\prod_{b=a+1}^{\la^{(l)}}\frac{[v^{(l)}_a-v^{(l)}_b-1]}{[v^{(l)}_a-v^{(l)}_b]}
\Biggr)
,
\en
where  we set $t^{(l)}_a=q^{2v^{(l)}_a }$, $t^{(N)}_s=z_s=q^{2u_s}$, $M=\sum_{l=1}^{N-1}\la^{(l)}$ 
and $C_{\mu_s,l+1}(s)=\sum_{j=s+1}^n\langle\bep_{\mu_j},h_{\mu_s,l+1}\rangle$. 
The symbol $<\ \ >^{Sym}$ denotes the symmetric part of the operator product expansion coefficient 
$<\ \ >$ defined by for $\cO=\Phi_N, F_l$, 
\be
&&\cO(w_a)\cO(w_b)=<\cO(w_a)\cO(w_b)>:\cO(w_a)\cO(w_b):. 
\en
 Hence $\widetilde{\Phi}(t,z)$ is an operator valued symmetric function in $z=(z_1,\cdots,z_n)$ as well as in $t^{(l)}=(t^{(l)}_1,\cdots,t^{(l)}_{\la^{(l)}})$ for each $l$. 
The function $\widetilde{W}_I(t,z,\Pi)$ is the elliptic weight function of type $\sln$. It is obtained by collecting all non-symmetric part of $<\ \ >$'s and factors arisen from the exchange among 
$\Phi_N$'s and  $F_l$'s and by making symmetrizations ${\rm Sym}_{t^{(l)}}$ of all entries in $ t^{(l)}$.
  Hence $\widetilde{W}_I(t,z,\Pi)$ is a symmetric function in $v^{(l)}_a \ (a=1,\cdots,\la^{(l)})$  for each $l$.  

The weight function $\widetilde{W}_I(t,z,\Pi)$ has several nice properties such as  
\begin{itemize}
\item{The triangular property}: 
For $I,J\in \cI_\la$, 
\begin{itemize}
\item[(1)] $\tW_{J}(z_I,z,\Pi)=0$ unless $I\leqslant J$.
\item[(2)] 
$\ds{\tW_{I}(z_I,z,\Pi)=\prod_{1\leq k<l\leq N}\prod_{a\in I_k}\prod_{b\in I_l\atop a<b}\frac{[u_b-u_a]}{[u_b-u_a+1]}}$
\end{itemize}
Here  $\leqslant$ denotes the partial ordering defined by 
\be
I\leqslant J \Leftrightarrow i^{(l)}_a \leq j^{(l)}_a\qquad \forall l, a. 
\en
for  $I^{(l)}=\{i^{(l)}_1<\cdots<i^{(l)}_{\la^{(l)}}\}$ and $J^{(l)}=\{j^{(l)}_1<\cdots<j^{(l)}_{\la^{(l)}}\}$ $(l=1,\cdots,N)$. 
We also denote by $t=z_I$ the specialization $t^{(l)}_a=z_{i^{(l)}_a}$ $(l=1,\cdots,N-1, a=1,\cdots,\la^{(l)})$\cite{RTV}.


\item {The transition property}: 
\be
&&
\widetilde{W}
_{I_{\cdots\ \mu_{i+1}\mu_{i}\cdots}}
(t, 
\cdots,z_{i+1},z_{i},\cdots
,\Pi)\nn\\
&&=
\sum_{\mu_{i}',\mu_{i+1}'}\bR(z_{i}/z_{i+1},\Pi q^{-2\sum_{j=i}^{n}\langle\bep_{\mu_j},h\rangle})_{\mu_{i}\mu_{i+1}}^{\mu_{i}'\mu_{i+1}'}\ 
\widetilde{W}_{I_{\cdots\ \mu'_i\mu'_{i+1}\cdots}}(t, 
\cdots,z_i,z_{i+1},\cdots
,\Pi)\lb{transsi}
\en
\end{itemize}  
as well as the orthogonality, quasi-periodicity, the shuffle product formula compatible to the wheel conditions etc.\cite{Konno17}. These properties are used to identify $\widetilde{W}_I(t,z,\Pi)$
 with the elliptic stable envelope for the cotangent bundle to the partial flag variety (Sec.\ref{secGR}).

\section{Integral Solution of the Elliptic $q$-KZ equation}\lb{secSolqKZ}
Let us consider the $n$-point function $\phi_{\mu_1\cdots\mu_n}(z_1,\cdots, z_n)$ with  the zero weight condition $\sum_{i=1}^{n}\bep_{\mu_i}=0$, which is equivalent to consider $\la=(m^N)$ with  $mN=n$. In this case, one can take a trace over the space $\F_{a,\nu}(\xi,\eta)$.  It then turns out that the trace gives a solution of the elliptic $q$-KZ equation due to the cyclic property of the trace and the commutation relation of the vertex operators \eqref{typeIcr}\cite{FJMMN,Konno17}. From \eqref{nptop}, one gets the following  elliptic hypergeometric integral. 
\be
&&\tr_{\F_{a,\nu}}\left(Q^{-d}\Phi_{\mu_1}(z_1)\cdots\Phi_{\mu_N}(z_N)\right)
=\oint_{\T^M
} \prod_{l=1}^{N-1}\prod_{a=1}^{
\la^{(l)}
}\frac{dt^{(l)}_a}{2\pi i t^{(l)}_a}
\ e(t,\Pi) {\Phi}(t,z)\tW_I(t,z,\Pi),\\
&&e(t,\Pi)=
\exp{\left\{\frac{\sum_{l=1}^{N-1}\log 
(
\Pi_l/\Pi_{l+1})
\sum_a\log 
(t^{(l)}_a)}{\log p}\right\}},\\
&& {\Phi}(t,z)=
\prod_{l=1}^{N-1}\left[
\prod_{a=1}^{
\la^{(l)}
}\prod_{b=1}^{
\la^{(l+1)}
}
\frac{\Gamma(t^{(l)}_a/t^{(l+1)}_b;p,Q)}{\Gamma(p^*t^{(l)}_a/t^{(l+1)}_b;p,Q)}
\prod_{1\leq a<b\leq \la^{(l)}}\frac{\Gamma(p^*t^{(l)}_a/t^{(l)}_b,p^*t^{(l)}_b/t^{(l)}_a;p,Q)}{\Gamma(t^{(l)}_a/t^{(l)}_b,t^{(l)}_b/t^{(l)}_a;p,Q)}\right]. 
\en
Note that only the weight function $\tW_I(t,z,\Pi)$ carries the tensor indices through $I$.  
We expect that $\tW_I(t,z,\Pi), {I\in \cI_{(m^N)}}$ form a basis of the twisted de Rham cohomology   associated with $T_{Q,z_i}\Phi(t,z)/\Phi(t,z)$, $i=1,\cdots,n$, where  $T_{Q,z_i}$ is the shift operator $T_{Q,z_i}f(\cdots,z_i,\cdots)=f(\cdots,Qz_i,\cdots)$. Furthermore, to specify the cycle of the integral one may insert to the integral an additional elliptic weight function $\tW_J(t,z,\Pi)$, $J\in \cI_{(m^N)}$, whose elliptic nome is $Q$ instead of $p$. Thus obtained elliptic hypergeometric integrals labeled by $I,J\in \cI_{(m^N)}$ are conjectured to give fundamental solutions of the elliptic $q$-KZ equation.  See \cite{TV97} for the trigonometric $\slt$ case. 
 
In the trigonometric limit $Q\to 0$, the trace goes to the following vacuum expectation value.
\be
&&\bra{0}\Phi_{\mu_1}(z_1)\cdots\Phi_{\mu_N}(z_N)\ket{0}
=\oint_{\T^M
} \prod_{l=1}^{N-1}\prod_{a=1}^{
\la^{(l)}
}\frac{dt^{(l)}_a}{2\pi i t^{(l)}_a}
\  e(t,\Pi){\Phi}^{trig.}(t,z)\tW_I(t,z,\Pi),\\
&& {\Phi}^{trig.}(t,z)=
\prod_{l=1}^{N-1}\left[
\prod_{a=1}^{
\la^{(l)}
}\prod_{b=1}^{
\la^{(l+1)}
}
\frac{(p^*t^{(l)}_a/t^{(l+1)}_b;p)_\infty}{(t^{(l)}_a/t^{(l+1)}_b;p)_\infty}
\prod_{1\leq a<b\leq 
\la^{(l)}
}\frac{(t^{(l)}_a/t^{(l)}_b,t^{(l)}_b/t^{(l)}_a;p)_\infty}{(p^*t^{(l)}_a/t^{(l)}_b,p^*t^{(l)}_b/t^{(l)}_a;p)_\infty}\right].
\en
This gives the integral formula for the vertex function for the equivariant $\mathrm{K}$-theory $\mathrm{K}_T(T^*fl_\la)$ with $T=(\C^\times)^n\times \C^\times$, which is the generating function of counting the quasi-maps  from $\PP^1$ to $T^*fl_\la$\cite{KPSZ}. Here the elliptic weight function $\tW_I$ plays the role of the pole subtraction matrix\cite{AO}.  
In particular, for the full flag variety case i.e.  $\la^{(l)}=l\ (l=1,\cdots,N-1)$ with $n=N$,  it gives the integral formula for the Macdonald symmetric function\cite{Mimachi}, which is also known as the 3d vortex partition function of the $\mathrm{T}[U(N)]$ theory, if one drops the weight function $\tW_I$ from the integrand.   Indeed the first and the second blocks in ${\Phi}(t,z)$ are elliptic analogues of the Cauchy type reproducing kernel and the Macdonald weight function, respectively, by reading $(p,p^*)$ as  $(q,t)$ in the Macdonald theory.   
Hence the vertex operators of $U_{q,p}(\slnh)$ play the role of  the vertex operators of the deformed $W$-algebra $\cW_{p,p^*}\sln$\cite{AS,FOS19}.
%
%

\section{Algebraic Analysis of the Elliptic Lattice Models}\lb{secAA}
The algebraic analysis \cite{JM} is a scheme of mathematical formulation of the  solvable 2d lattice models defined by the $R$-matrices.  It formulates the model on the infinite lattice directly by using  both finite and infinite dimensional representations of the corresponding affine quantum groups $U_q(\g)$ or  elliptic quantum groups $U_{q,p}(\g)$. 
In particular, the infinite dimensional representations are identified 
with the spaces of states of the model specified by  the ground state boundary condition, where  the configurations at  enough far sites from the center of the  lattice is given by  the ground state configurations. 
In addition,  the two types intertwining operators, the type I and II vertex operators, of the  quantum groups  realize the local operators of the model such as spin operators  
acting on the infinite dimensional representations (type I)  and   the creation operators of the physical excitations (type II), respectively.    
Then the traces of compositions of the vertex operators give correlation functions as well as form factors of the models. Hence these are characterized as solutions of the (elliptic) $q$-KZ equations\cite{IFrRe,FJMMN,JM} as discussed in the last section for the type I case. In this sense, the algebraic analysis is an off-critical extension of CFT
\cite{BPZ}.  

There are series of elliptic  solvable lattice models defined by the elliptic solutions of the face type YBE, i.e. the elliptic dynamical $R$-matrices, associated with affine Lie algebras\cite{ABF,JMO,DJKMO,Kuniba,KS}.  
The critical behavior of such models  is  described by the 
$W$-algebra $\cW\bar{\g}$ of the coset type $(\g)_{r-h^\vee-k}\oplus (\g)_{k}\supset (\g_{\rm diag.})_{r-h^\vee}$\cite{DJO,DJKMO}. 
Here $r$ is taken as a  real positive parameter and called the restriction height. 
Accordingly the elliptic quantum group $U_{q,p}(\g)$  plays a central role in the algebraic analysis of these models as a deformation  of the $W$-algebra $\Wppg
$. See Sec.\ref{UqpWpps}. The spaces of states with  the ground state boundary condition are identified with the irreducible $\Wppg
$-modules, whose direct sum  are isomorphic to the $U_{q,p}(\g)$-modules (Sec.4.2). The type I and II vertex operators of $U_{q,p}(\g)$ provide correlation functions and form factors of the model. See \cite{FJMMN,LP,KKW}. 
  
 In the case $\g=\slnh$,  one also has the vertex type elliptic lattice models\cite{Baxter,Belavin}. These    are related to the face models by the so-called vertex-face correspondence\cite{Baxter,JMO87}. The algebraic analysis of the vertex models is carried out based on the one of the face models by applying the vertex-face correspondence\cite{LasPu,KKW,CKSW2}.  
%
%

\section{Level-0 Action on the Gelfand-Tsetlin Basis }\lb{secGTbasis}
Let $(\pi_z,\hV_{z})$ be the evaluation representation of $U_{q,p}(\slnh)$ in Sec.\ref{vecrep}. 
The action of the $L$-operator $L^+(1/w)$ on $\hV_w\tot\hV_{z}$ is extended to the tensor space $ \hV_w\tot  \hV_{z_1} \tot \cdots \tot \hV_{z_n}$ by the comultiplication 
\be
&&(\pi_{z_1}\otimes \cdots\otimes \pi_{z_n}){\Delta'}^{(n-1)}(L^\pm(1/w)) \\
&&=\bar{R}^{(0n)}(z_n/w,\Pi^* q^{2\sum_{j=1}^{n-1}h^{(j)}}){\bR}^{(0n-1)}(z_{n-1}/w,\Pi^* q^{2\sum_{j=1}^{n-2}h^{(j)}})\cdots \bar{R}^{(01)}(z_1/w,\Pi^*).
\en 
Here we take the opposite comultiplication $\Delta'$ following \cite{Konno18ESE}.   
This is the action on the standard basis  
\be
&&v_I
:=v_{\mu_1}\tot \cdots \tot v_{\mu_n},\qquad I=I_{\mu_1\cdots \mu_n}\in \cI_{\la},\ \la\in \N^N,\ |\la|=n.
\en

Let us consider the following change of basis\cite{Konno18ESE}. 
\bea
&&\xi_I=\sum_{J\in \cI_{\la}}\tW_J(z^{-1}_I,z^{-1},\Pi q^{2\sum_{j=1}^n<\bep_{\mu_j},h>})
v_J.\lb{xiv}
\ena
The transition matrix is given by the specialized elliptic weight function so that it is lower triangular (Sec.\ref{secWF}).   It turns out that the new vectors $\xi_I$ 
$(I\in \cI_\la, \la\in \N^N, |\la|=n)$ form the Gelfand-Tsetlin (GT) basis of $\hV_{z_1}\tot \cdots \tot \hV_{z_n}$,  
on which the commuting subalgebra 
generated by $\al_{j,m}$ $(j=1,\cdots,N-1, \ m\in \Z\backslash\{0\} )$ is diagonalized simultaneously.  
In fact the  level-0 action of the elliptic currents of $U_{q,p}(\slnh)$ on 
$\xi_I$ is given by 
\be
&&\phi_j^\pm(w)\xi_I=
\left.\prod_{a\in I_j}\frac{[u_a-v+1]}{[u_a-v]}\right|_{\pm}\left.\prod_{b\in {I_{j+1}}}\frac{[u_b-v-1]}{[u_b-v]}\right|_{\pm}\ e^{-Q_{\al_j}}\xi_I\\
&&E_{j}(w)\xi_I=a^*
\sum_{i\in I_{j+1}}\delta(z_i/w)\prod_{k\in I_{j+1}\atop \not=i}\frac{[u_i-u_k+1]}{[u_i-u_k]}\ e^{-Q_{\al_j}}\xi_{I^{i'}}\\
&&F_{j}(w)\xi_I=a\sum_{i\in I_{j}}\delta(z_i/w)\prod_{k\in I_{j}\atop \not=i}\frac{[u_k-u_i+1]}{[u_k-u_i]}\ \xi_{I^{'i}}\\[-5mm]
\en
where ${aa^*=-\frac{1}{q-q^{-1} }\frac{[0]'}{[1]}}$, and $I=(I_1,\cdots, I_N)$, $I^{i'}\in \cI_{(\la_1,\cdots,\la_j+1,\la_{j+1}-1,\cdots,\la_N)}$, \\
$I^{'i}\in \cI_{(\la_1,\cdots,\la_j-1,\la_{j+1}+1,\cdots,\la_N)}$  with 
\be
&&(I^{i'})_j=I_j\cup \{i\},\quad (I^{i'})_{j+1}=I_{j+1}-\{i\},\quad (I^{i'})_k=I_k \ \ (k\not=j,j+1), \\
&&(I^{'i})_j=I_j- \{i\},\quad (I^{'i})_{j+1}=I_{j+1}\cup\{i\},\quad (I^{'i})_k=I_k \ \ (k\not=j,j+1).
\en
Note that this action is given in terms of the partitions in completely a combinatorial way.  In 
the next section, this is identified with a geometric action of $U_{q,p}(\slnh)$ on the 
equivariant elliptic cohomology $\ds{\bigoplus_{\la\in \N^N,|\la|=n}\E_T(T^*fl_\la)}$ under the 
identification of the elliptic weight function with the elliptic stable envelope.

\section{Geometric Interpritation}\lb{secGR}
The notion of stable envelopes was initiated in \cite{MO,Okounkov}. They form a good basis of  equivariant cohomology and  $\mathrm{K}$-theory of  Nakajima quiver varieties\cite{Na94,Na98} and has nice applications to enumerative geometry, geometric representation theory and quantum integrable systems. In particular, the transition matrices between stable envelopes defined for different chambers give the $R$-matrices satisfying the YBE so that  they provide a new geometric formulation of quantum groups as well as quantum integrable systems associated with the quiver varieties.  The elliptic version of stable envelopes  was defined in \cite{AO} for the equivariant elliptic cohomology $\E_T(X)$ of the Nakajima quiver varieties $X$. Remarkably 
they provide the elliptic dynamical $R$-matrix as their transition matrix, which coincides  for example with \eqref{ellR} for the case $X=T^*fl_\la$. We here expose an identification of the elliptic stable envelopes for $\E_T(T^*fl_\la)$ with the elliptic weight functions for $U_{q,p}(\slnh)$ described in Sec.\ref{secWF}. This makes more transparent the correspondence between geometry of quiver variety and 
 representation of EQG. For the definition of the elliptic cohomology and the elliptic stable envelopes, the reader can consult\cite{AO}.

The elliptic stable envelope $\Stab_{\gC}(F_I)$ for $\E_T(T^*fl_\la)$ 
at the torus fixed point $F_I$ labeled by the partition $I\in \cI_\la$ can be constructed  explicitly by the abelianization procedure\cite{AO,Konno18ESE}. 
On the other hand, let us consider the elliptic weight function modified as  
\be
&&\cW_I(t,z,\Pi)=\frac{H_\la(t,z)\tW_I(t,z,\Pi)}{E_\la(t)},\\
&&H_\la(t,z)=\prod_{l=1}^{N-1}\prod_{a=1}^{\la^{(l)}}\prod_{b=1}^{\la^{(l+1)}}\left[v^{(l+1)}_b-v^{(l)}_a+1\right],\lb{Efunc}
\quad E_\la(t)=\prod_{l=1}^{N-1}\prod_{a=1}^{\la^{(l)}}\prod_{b=1}^{\la^{(l)}}[v^{(l)}_b-v^{(l)}_a+1].\lb{def:Ela}
\en
This yields the expresion 
\be
\cW_I(t,z,\Pi)
&=&{\rm Sym}_{t^{(1)}}\cdots {\rm Sym}_{t^{(N-1)}}\ {U}_I(t,z,\Pi),\lb{def:cW}\\
U_I(t,z,\Pi)&=&\prod_{l=1}^{N-1}\frac{\prod_{a=1}^{\la^{(l)}} 
    u^{(l)}_I(t^{(l)}_a,t^{(l+1)},\Pi_{\mu_{\mbox{\tiny${i^{(l)}_a}$}},l+1}q^{-2C_{\mu_{\mbox{\tiny${i^{(l)}_a}$}},l+1}(i^{(l)}_a)})}{\prod_{1\leq a<b\leq \la^{(l)}}{[v^{(l)}_a-v^{(l)}_b]}{[v^{(l)}_b-v^{(l)}_a-1]}},\lb{def:cU}\\
u^{(l)}_I(t^{(l)}_a,t^{(l+1)},\Pi_{j,k})&=& 
\left.
\frac{\left[v^{(l+1)}_b-v^{(l)}_a+(P+h)_{j,k}
\right]}{[(P+h)_{j,k}
]}
\right|_{i^{(l+1)}_b=i^{(l)}_a}
\nn\\
&&\qquad\times 
\prod_{b=1\atop i^{(l+1)}_b>i^{(l)}_a}^{\la^{(l+1)}}{\left[v^{(l+1)}_b-v^{(l)}_a\right]}
\prod_{b=1\atop i^{(l+1)}_b<i^{(l)}_a}^{\la^{(l+1)}}{\left[v^{(l+1)}_b-v^{(l)}_a+{1}\right]}.
\lb{def:u}
\en


Comparing these expressions, one finds $\Stab_{\gC}(F_I)$ and $\cW_I(z,t,\Pi)$ are identical as
\bea
&&\Stab_{\gC}(F_I)= \cW_{\sigma_0(I)}(\tit,\sigma_0(z^{-1}),\Pi^{-1}). 
\ena
Here for the longest element $\sigma_0$ in $\gS_n$ and $I=I_{\mu_1\cdots\mu_n}$, we define $\s_0(I)=I_{\mu_{\s_0(1)}\cdots\mu_{\s_0(n)}}$ and set $\s_0(I)=(\tilde{I}_1,\dots,\tilde{I}_N)$, $\tilde{I}^{(l)}=\tilde{I}_1\cup \cdots \cup \tilde{I}_l$  
and $\tilde{I}^{(l)}=\{\tilde{i}^{(l)}_1<\cdots<\tilde{i}^{(l)}_{\la^{(l)}}\}$.  Then  $\tit$ denotes the set of all   
$t^{(l)}_a$ corresponding to $\tilde{i}^{(l)}_a$. 
 The chamber $\gC$ is taken as  $|z_1|<\cdots<|z_n|$. 
In addition, the restriction $t=z_J$ defined in Sec.\ref{secWF} can be identified with the restriction of the stable class $\Stab_{\gC}(F_I)$ to  
 the fixed point $F_J$. One then finds  
\bea
&&{\rm Stab}_\gC(F_I)\vert_{F_J}=\cW_{\sigma_0(I)}(z^{-1}_{J},\sigma_{0}(z^{-1}), \Pi^{-1}).\lb{StabWFP}
\ena
Such coincidence as well as the integral formula for the vertex function
 obtained as correlation function of the  vertex operators
   in Sec.\ref{secSolqKZ} lead us to the following dictionary from 
the level-0 representation of $U_{q,p}(\slnh)$
 to the equivariant elliptic cohomology $\E_T(T^*fl_\la)$. 
\be
 \mbox{evaluation parameters}\ z_1,\cdots,z_n \ &\leftrightarrow& \ \mbox{equivariant parameters}\ 
 \in (\C^\times)^n\\[1mm]
 \mbox{dynamical parameters}\ \Pi_j  \ &\leftrightarrow& \ \mbox{K\"ahler parameters}\in \mathrm{Pic}_T(X)\otimes_\Z E \qquad\\[1mm]
  \mbox{integration variables} \ &\leftrightarrow& \ \mbox{Chern roots of the tautological vector bundle $\cV_l$}\ 
\\
\ t^{(l)}_a\ (a=1,\cdots,\la^{(l)}) \quad &
& \ \mbox{ of $\rk\; \la^{(l)}$} \mbox{on the $l$-th vertex of the quiver diagram}
\\[1mm]
\mbox{the weight}\ (\mu_1,\cdots,\mu_n) \ &\leftrightarrow& \ \mbox{the fixed point}\ F_{I_{\mu_1,\cdots,\mu_n} } \ 
\en

Furthermore, denoting the fixed point classes by $[I]$  one has the localization formula for the stable envelope. 
\be
&&{\rm Stab}_\gC(F_K)=\sum_{I\in \cI_\la}\frac{{\rm Stab}_\gC(F_K)\vert_{F_I}}{N(z_I)}\ [I],\lb{StabI}
\en
where  $N(z_I)$ is a certain normalization. Under the identification \eqref{StabWFP}, one can compare this with the change of basis relation from the GT basis $\xi_I$ to the standard basis $v_K$, i.e. the reversed relation of \eqref{xiv} obtained by appying the orthogonality property of the elliptic weight functions.  Then one finds that these two are identical under the correspondence\cite{Konno18ESE}
 \be
  \mbox{standard bases}\ v_I
     \ &\leftrightarrow& \ \mbox{stable classes }\  \Stab_{\gC}(F_I) \\[1mm]
 \mbox{the Gelfand-Tsetlin bases} \ \xi_I\ &\leftrightarrow& \ \mbox{fixed point classes}\ [I].
\en
This allows us to interpret the level-0 action of $U_{q,p}(\slnh)$ on the GT basis of $\cV_{z_1}\otimes \cdots \otimes \cV_{z_n}$  in Sec.\ref{secGTbasis} as an action on the fixed point classes in $\bigoplus_{\la\in\N^N,|\la|=n}\E_T(T^*fl_\la)$.  



\section{Elliptic Quantum Toroidal Algebras}\lb{secEQTA}
A toroidal algebra $\g_{tor}$ associated with a complex simple Lie algebra $\bar{\g}$ 
 is a canonical two-dimensional central extension of the double loop Lie algebra $\C^{\times}\times \C^{\times} \to \bar{\g}$. Its quantum version, i.e. quantum toroidal algebra, was introduced in \cite{GKV}. 
 We denote it by $U_{q,\kappa}(\g_{tor})$.  
It is formulated as an extension of 
the quantum affine algebra $U_q({\g})$ in  the Drinfeld  realization\cite{DrinfeldNew} by 
introducing the new generators, i.e. the 0-th generators w.r.t. the index set $I$,  and  replacing the finite type Cartan matrix $(a_{i,j})\ i,j\in \bI$ appearing in the defining relation of $U_q({\g})$ with the affine type generalized Cartan matrix $(a_{i,j})\ i,j\in I$. This procedure is called the quantum affinization. 
 In particular, for the case $\bar{\g}=\gl_N$ 
one can introduce an extra parameter $\kappa$ due to a cyclic property of the affine Dynkin diagram of  type $A^{(1)}_{N-1}$. 

In the same way, the elliptic quantum toroidal algebra $U_{q,\kappa,p}(\g_{tor})$ is formulated as the quantum affinization of the elliptic algebra $U_{q,p}(\g)$ by adding the 
0-th generators $\al_{0,m}, e_{0,n}, f_{0,n}, K^\pm_0\ (m\in\Z\backslash\{0\}, n\in \Z)$\cite{KO23}. As a result, the element $\prod_{i\in I}(K^+_i)^{a_i^\vee}$ becomes a central element. Here 
 $a^\vee_i\ (i\in I)
 $ denote the colabel of the affine Dynkin diagram\cite{Kac}. 
 In addition,  ${U}_{q,\kappa,p}(\g_{tor})$ possesses two  subalgebras,  both of which are isomorphic to the elliptic  algebra $U_{q,p}(\g)$: 
the one  generated by $\al_{i,l}, k_i^\pm, e_{i,n}, f_{i,n}, q^{\pm {c}/{2}}, 
q^d$ $(i\in \bI,  l\in \Z\backslash \{0\}, n\in \Z  
)$ and the other by $e_{i,0}, f_{i,0},  K^{\pm}_i, q^d$ $(i\in I)$. 
These are elliptic analogues of those in  $U_{q,\kappa}(\g_{tor})$\cite{GKV}.

On the other hand, the $\gl_1$ version of the quantum toroidal algebra $U_{q,t}(\gl_{1,tor})$  was introduced in \cite{Miki4} as a $q,t$-deformation of the $W_{1+\infty}$ algebra.
Remarkably, $U_{q,t}(\gl_{1,tor})$ is  isomorphic to the elliptic Hall algebra\cite{Sc1}.  It is also remarkable that   
representations of $U_{q,t}(\gl_{1,tor})$ have a deep connection to the Macdonald theory\cite{Miki4,FT, 
SV1,FOS19}. 
 In addition, representations of  $U_{q,t}(\gl_{1,tor})$ provide a relevant scheme of  
  calculation of the instanton partition functions\cite{Nekrasov04} as well as  a study of  
 the Alday-Gaiotto-Tachikawa (AGT) correspondence\cite{AGT} for the 5d and 6d lifts of the 4d ${\cN=2}$ SUSY gauge theories, which are known as  linear quiver gauge theories.

The elliptic  quantum toroidal algebra  $U_{q,t,p}(\gl_{1,tor})$ was formulated in \cite{KOgl1}. It has several nice representations such as  
the level-(0,0) representation realized in terms of the elliptic Ruijsenaars difference operators\cite{Ruijs} as well as an elliptic analogue of the $q$-Fock representation of $U_{q,t}(\gl_{1,tor})$ 
which is expected to give a geometric representation on the equivariant elliptic cohomology of the Hilbert scheme of points on $\C^2$, $\bigoplus_{n}\E_T(\mathrm{Hilb}_{n}(\C^2))$. 
These  suggest a deep connection of $U_{q,t,p}(\gl_{1,tor})$ 
to an expected elliptic analogue of the theory of Macdonald symmetric function. 
 It is also remarkable that $U_{q,t,p}(\gl_{1,tor})$ realizes a deformation of the Jordan quiver $W$-algebra $\cW_{p,p^*}(\Gamma(\widehat{A}_0))$ introduced in  \cite{KimPes} in the same way as the elliptic algebra $U_{q,p}(\g)$ realizes the deformed $W$-algebra $\Wppg
$ discussed in Sec.\ref{UqpWpps}. Namely, the level-$(1,N)$ elliptic currents of $U_{q,t,p}(\gl_{1,tor})$ 
gives a realization of the screening currents and a composition of the vertex operators 
 w.r.t. the Drinfeld comultiplication gives a realization of the generating functions of $\cW_{p,p^*}(\Gamma(\widehat{A}_0))$. Moreover thus obtained realization  turns out to provide a  relevant scheme to the  instanton calculus of  
 the 5d and 6d lifts of the  4d $\cN=2^*$ SUSY gauge theories, i.e.  the 
  $\cN=2$ SUSY gauge theories coupled with the adjoint matter\cite{Nekrasov04}, known as 
 the Jordan quiver gauge theories\footnote{It is also called the ADHM quiver gauge theories.}.  There one of the essential feature is that $U_{q,t,p}(\gl_{1,tor})$ at the level (1,$N$) possesses  
 the  four parameters $q, t, p, p^*$ satisfying one constraint $p/p^*=t/q$,  
which play  the role of  the $SU(4)$ $\Omega$-deformation parameters\cite{Nekrasov}.

\section{Summary}
We have exposed the EQGs $U_{q,p}(\g), U_{q,\kappa,p}(\g_{tor})$ and $U_{q,t,p}(\gl_{1,tor})$ in the Drinfeld realization. The main features of them are 
summarized as follows.
\begin{itemize} 
%
\item They allow to treat both the level-0 and level $\not=0$ representations in a unified way. In particular the level $\not=0$ representation yields the two elliptic nomes $p, p^*$, 
which are identified with $q,t$ in the Macdonald theory and the deformed $W$-algebras. In the same way $U_{q,t,p}(\gl_{1,to})$ possesses the $SU(4)$ $\Omega$-deformation parameters $q,t,p,p^*$ in the level $(1,N)$ representation. 

\item The coalgebra structure of them allows to formulate the vertex operators as intertwining operators of modules, where the level-0 and level $\not=0$ representations are mixed as tensor product. 

\item The vertex operators are important objects having many applications to 
formulations of mathematical and physical quantities such as elliptic weight functions, integral solutions of the elliptic $q$-KZ equation, vertex functions, correlation functions and form factors in algebraic analysis of elliptic solvable lattice models etc.. 
\item They realize 
a deformation of the W-algebras including the affine quiver types. 
In particular, the vertex operators of EQGs provide those of the deformed $W$-algebras. 
\end{itemize}

We also have shown a correspondence between the level-0 representation of $U_{q,p}(\slnh)$ and the equivariant elliptic cohomology $\E_T(T^*fl_\la)$ by 
deriving the elliptic weight functions, the integral solutions of the (elliptic) $q$-KZ equation as well as the level-0 action on the Gelfand-Tsetlin basis.
In this process, the vertex operators of  $U_{q,p}(\slnh)$ play an essential role. 
However, in the geometry side neither the non-zero level representations nor the vertex operators intertwining them have yet been understood.  It should be outstanding to formulate the vertex operators of EQGs geometrically.  
  
In addition, knowing that the formulations of elliptic quantum groups and equivariant elliptic cohomologies possess the equivariant  and K\" ahler parameters in a unified way, one of their promising application seems to be the 3d-mirror symmetry for  a  pair of  supersymmetric gauge theories with $\cN$ = 4 supersymmetry\cite{AO}. This provides an important example of the deeper mathematical problem of symplectic duality for conical symplectic resolutions. 
In such gauge theories, Nakajima quiver varieties are known to give moduli spaces of vacua called the Higgs branch. The  3d-mirror symmetry  states   
that the Higgs branch of the dual theory 
conjecturally coincides with the Coulomb branch of the original one in such a way that 
the roles of the equivariant  and the K\"ahler parameters of the dual theories are exchanged. 
To show the duality at the level of physical quantities such as vertex functions would be outstanding.


\end{document}